\documentclass[reqno]{amsart}%
\usepackage{float}
\usepackage{amsmath}
\usepackage{graphicx}
\usepackage{latexsym}
\usepackage{amsfonts}
\usepackage{amssymb}%
\newtheorem{theorem}{Theorem}

\newtheorem{lemma}[theorem]{Lemma}

\theoremstyle{definition}

\newtheorem{definition}[theorem]{Definition}
\newtheorem{example}[theorem]{Example}

\newtheorem{remark}[theorem]{Remark}

\DeclareSymbolFont{letters}{OML}{cmm}{m}{it}
\DeclareSymbolFont{bletters}{OML}{cmm}{b}{it}
\DeclareMathSymbol{\bphi}{\mathord}{bletters}{"1E}
\DeclareMathSymbol{\bvarphi}{\mathord}{bletters}{"27}
\DeclareMathSymbol{\bvartheta}{\mathord}{bletters}{"23}
\DeclareMathSymbol{\bvareps}{\mathord}{bletters}{"22}
\DeclareMathSymbol{\bmu}{\mathord}{bletters}{"16}
\DeclareMathSymbol{\bpsi}{\mathord}{bletters}{"20}
\DeclareMathSymbol{\bpi}{\mathord}{bletters}{"19}
\DeclareMathSymbol{\bnu}{\mathord}{bletters}{"17}
\DeclareMathSymbol{\balpha}{\mathord}{bletters}{"0B}
\DeclareMathSymbol{\bbeta}{\mathord}{bletters}{"0C}
\DeclareMathSymbol{\bgamma}{\mathord}{bletters}{"0D}
\DeclareMathSymbol{\blambda}{\mathord}{bletters}{"15}
\DeclareMathSymbol{\bomega}{\mathord}{bletters}{"21}
\DeclareMathSymbol{\bkappa}{\mathord}{bletters}{"14}
\DeclareMathSymbol{\bLambda}{\mathord}{bletters}{"03}
\DeclareMathSymbol{\bsigma}{\mathord}{bletters}{"1B}
\DeclareMathSymbol{\btheta}{\mathord}{bletters}{"12}
\DeclareMathSymbol{\bTheta}{\mathalpha}{operators}{"02}
\DeclareMathSymbol{\bfeta}{\mathord}{bletters}{"11}
\DeclareMathSymbol{\bzeta}{\mathord}{bletters}{"10}
\begin{document}
\title[Lower deviations for supercritical GWP]{Lower deviation probabilities for supercritical Galton-Watson processes}
\author[Fleischmann]{Klaus Fleischmann}
\address{Weierstrass Institute for Applied Analysis and Stochastics, Mohrenstr.\ 39,
D--10117 Berlin, Germany}
\email{fleischm@wias-berlin.de}
\urladdr{http://www.wias-berlin.de/\symbol{126}fleischm}
\author[Wachtel]{Vitali Wachtel}
\address{Weierstrass Institute for Applied Analysis and Stochastics, Mohrenstr.\ 39,
D--10117 Berlin, Germany}
\email{vakhtel@wias-berlin.de}
\urladdr{http://www.wias-berlin.de/\symbol{126}vakhtel}
\thanks{Supported by the DFG}
\thanks{Corresponding author: K.\ Fleischmann}
\keywords{Supercritical Galton-Watson process, local limit theorem, large deviation,
Cram\'er transform, concentration function, Schr\"{o}der equation,
B\"{o}ttcher equation}
\subjclass{Primary 60\thinspace J\thinspace80; Secondary 60\thinspace F\thinspace10.}
\date{April 28, 2005;\quad WIAS Preprint No. 1025;\quad LowerDev20.tex\quad ISSN
0946 -- 8633}
\maketitle

\thispagestyle{empty}\setcounter{page}{0}

\vfill\newpage

$\quad$\vfill

\begin{quotation}
\textsc{Abstract. }{\small There is a well-known sequence of constants}
$\,c_{n}$\thinspace\ describing the growth of supercritical Galton-Watson
processes $\,Z_{n\,}.$\thinspace\ With \textquotedblleft lower deviation
probabilities\textquotedblright\ we refer to $\,\mathbf{P}(Z_{n}=k_{n}%
)$\thinspace\ with $\,k_{n}=o(c_{n})$\thinspace\ as $\,n$\thinspace
\ increases. We give a detailed picture of the asymptotic behavior of such
lower deviation probabilities. This complements and corrects results known
from the literature concerning special cases. Knowledge on lower deviation
probabilities is needed to describe large deviations of the ratio
$\,Z_{n+1}/Z_{n\,}.$\thinspace\ The latter are important in statistical
inference to estimate the offspring mean. For our proofs, we adapt the
well-known Cram\'{e}r method for proving large deviations of sums of
independent variables to our needs.
\end{quotation}

\begin{quote}
\vfill

{\scriptsize
\tableofcontents
}

\vfill\newpage
\end{quote}

\section{Introduction and statement of results}

\subsection{On the growth of supercritical processes\label{SS.growth}}

Let $Z=(Z_{n})_{n\geq0}$ denote a Galton-Watson process with offspring
generating function%
\begin{equation}
f(s)\ =\ \sum_{j\geq0}p_{j}s^{j},\qquad0\leq s\leq1, \label{not.f.p}%
\end{equation}
which is required to be non-degenerate, that is, $p_{j}<1,$\thinspace
\ $j\geq0.$\thinspace\ Suppose that $Z$ is supercritical, i.e. $f^{\prime
}(1)=:m\in(1,\infty).$ For simplicity, the initial state $Z_{0}\geq1$ is
always assumed to be deterministic, and, if not noted otherwise (as by an
application of the Markov property), we set $Z_{0}=1.$

It is well-known (see, e.g., Asmussen and Hering (1983) \cite[\S \thinspace
3.5]{AsmussenHering1983}) that
\begin{equation}
\text{there are }\,c_{n}>0\,\ \text{such that a.s.}\,\ c_{n}^{-1}%
Z_{n}\;\underset{n\uparrow\infty}{\longrightarrow}\;\text{some non-degenerate
}W. \label{not.cn}%
\end{equation}
In this sense, the sequence of constants $c_{n}$ describes the order of growth
of $\,Z.$\thinspace\ But, $\mathbf{P}(W=0)=q,$ with $q\in\lbrack0,1)$ the
smallest root of $\,f(s)=s,$ that is, the extinction probability of
\thinspace$Z$. On the other hand, $W$ restricted to $(0,\infty)\ $has a
(strictly) positive continuous density function denoted by $\,w.$%
\thinspace\ Therefore the following \emph{global limit theorem}\/ holds:%
\begin{equation}
\lim_{n\uparrow\infty}\mathbf{P}(Z_{n}\geq xc_{n})\ =\ \int_{x}^{\infty
}w(t)\,\mathrm{d}t,\qquad x>0. \label{ILT}%
\end{equation}
The normalizing sequence $(c_{n})_{n\geq0}$ can be chosen to have the
following additional properties:
\begin{subequations}
\label{properties.cn}%
\begin{align}
&  c_{0}=1\,\ \text{and}\,\ c_{n}\,<\,c_{n+1}\,\leq\,m\,c_{n\,},\,\ n\geq
0,\label{cn1}\\
&  c_{n}\,=\,m^{n}\,L(m^{n})\,\ \text{with }L\text{ slowly varying at
infinity,}\label{cn2}\\
&  \lim_{x\uparrow\infty}L(x)\,\ \text{exists; it is positive if and only
if}\,\ \mathbf{E}Z_{1\!}\log Z_{1}<\infty. \label{cn3}%
\end{align}
Because of (\ref{properties.cn}b,c), we may (and subsequently shall) take
\end{subequations}
\begin{equation}
c_{n}\,:=\,m^{n}\quad\text{if}\quad\mathbf{E}Z_{1\!}\log Z_{1}<\infty
.\vspace{2pt} \label{not.spec.cn}%
\end{equation}

\subsection{Asymptotic local behavior of $Z,$ purpose}

A local limit theorem related to (\ref{ILT}) is due to Dubuc and Seneta (1976)
\cite{DubucSeneta1976}, see also \cite[\S 3.7]{AsmussenHering1983}. To state
it we need the following definition.

\begin{definition}
[\textbf{Type }$(d,\mu)$]\label{D.type}We say the offspring generating
function $f$\thinspace\ \emph{is of type} $(d,\mu)$, if $\,d\geq1$%
\thinspace\ is the greatest common divisor of the set $\,\{j-\ell:\,j\neq
\ell,\ p_{j}p_{\ell}>0\},$\thinspace\ and $\,\mu\geq0$\thinspace\ is the
minimal \thinspace$j$\thinspace\ for which \thinspace$p_{j}>0$.\hfill
$\Diamond$
\end{definition}

Here is the announced \emph{local limit theorem}. Suppose $f$\thinspace\ is of
type $(d,\mu)$. Take $\,x>0,$\thinspace\ and consider integers\/ $\,k_{n}%
\geq1$ such that\/ $k_{n}/c_{n}\rightarrow x$\thinspace\ as $\,n\uparrow
\infty.$\thinspace\ Then, for each $\,j\geq1,$%
\begin{equation}
\lim_{n\uparrow\infty}\Bigl(c_{n}\,\mathbf{P}\!\left\{  Z_{n}=k_{n}%
\,\big|\,Z_{0}=j\right\}  -\,d\,\mathsf{1}_{\left\{  k_{n}\equiv j\mu
^{n}(\mathrm{mod}\,d)\right\}  }\,w_{j}(x)\Bigr)\ =\ 0, \label{A}%
\end{equation}
where\/ $\,w_{j}:=\sum\limits_{\ell=1}^{j}\left(
\genfrac{}{}{0pt}{1}{j}{\ell}%
\right)  q^{j-\ell}\,w^{\ast\ell}.$\bigskip

\noindent In particular, in our standard case $\,Z_{0}=1$ and if additionally
$\,k_{n}\equiv\mu^{n}\ (\mathrm{mod}\,d),$\thinspace\ then
\begin{equation}
\mathbf{P}(Z_{n}=k_{n})\ \sim\ d\,c_{n}^{-1}\,w(k_{n}/c_{n})\quad
\text{as}\,\ n\uparrow\infty\label{Intuition}%
\end{equation}
(with the usual meaning of the symbol $\sim$ as the ratio converges to $1).$

Statement (\ref{A}) [and especially (\ref{Intuition})] can be considered as
describing the local behavior of supercritical Galton-Watson processes in the
region of \emph{normal}\/ deviations (from the growth of the $\,c_{n\,}%
;$\thinspace\ `deviations' are meant here in a multiplicative sense, related
to the multiplicative nature of branching). But what about $\,\mathbf{P}%
(Z_{n}=k_{n})$\thinspace\ when $\,k_{n}/c_{n}\rightarrow0$ or $\infty
\,?$\thinspace\ In these cases we speak of \emph{lower}\/ and \emph{upper}\/
(local) deviation probabilities, respectively.

Lower deviations of $\,Z_{n}$\thinspace\ are closely related to large
deviations of $\,Z_{n+1}/Z_{n}$\thinspace\ (see Ney and Vidyashankar (2004)
\cite[Section~2.3]{NeyVidyashankar2004}). The latter are important in
statistical inference for supercritical Galton-Watson processes, since
$\,Z_{n+1}/Z_{n}$\thinspace\ is the well-known Lotka-Nagaev estimator of the
offspring mean.

The \emph{main purpose}\/ of the present paper is to study lower deviation
probabilities in their own and to provide a detailed picture (see
Theorems~\ref{T.Schroeder} and~\ref{T.Boettcher} below). As a starting point
we discuss a relevant claim in \cite{NeyVidyashankar2004} concerning an
important special case (see Sections~\ref{SS.lit} and~\ref{SS.contr} below).
Applications of our results for large deviations of $\,Z_{n+1}/Z_{n}%
$\thinspace\ and also to subcritical Galton-Watson processes are postponed to
a future paper.

Here is the program for the remaining introduction. After introducing a basic
dichotomy, we review in Sections~\ref{SS.lit} and~\ref{SS.contr} what is known
on lower deviations from the literature, before we state our results in
Sections~\ref{SS.Schr} and~\ref{SS.Boe}.

\subsection{A dichotomy for supercritical processes\label{SS.dich}}

Recalling that $f$ denotes the offspring generating function, $q$ the
extinction probability, and $m$ the mean,%
\begin{equation}
\text{set}\quad\gamma\,:=\,f^{\prime}(q),\quad\text{and define}\,\ \alpha
\,\ \text{by}\quad\gamma\,=\,m^{-\alpha}. \label{not.gamma.alpha}%
\end{equation}
Note that $\,\gamma\in\lbrack0,1)$\thinspace\ and $\,\alpha\in(0,\infty
].$\thinspace\ We introduce the following notion, reflecting a crucial
dichotomy for supercritical Galton-\hspace{-1.6pt}Watson processes.

\begin{definition}
[\textbf{Schr\"{o}der and B\"{o}ttcher case}]\label{D.dichot}\hspace{-2pt}For
our supercritical offspring law we distinguish between the \emph{Schr\"{o}der}%
\/ and the \emph{B\"{o}ttcher}\/ case, in dependence on whether\/
$\,p_{0}+p_{1}>0$\thinspace\ or $\,=0.$\hfill$\Diamond$
\end{definition}

\noindent Obviously, $f$\thinspace\ is of Schr\"{o}der type if and only if
$\,\gamma>0,$\thinspace\ if and only if $\,\alpha<\infty.$

Next we want to collect a few basic facts from the literature concerning that
dichotomy. Clearly, $f$ can be considered as a function on $D,$ where $D$
denotes the closed unit disc in the complex plane. As usual, denote by
$\,f_{n}$\thinspace\ the $n^{\mathrm{th}}$ iterate of $\,f.$\thinspace\ We
start with the \emph{Schr\"{o}der case}. Here it is well-known (see,
e.g.,\ \cite[Lemma~3.7.2 and Corollary~3.7.3]{AsmussenHering1983}) that
\begin{equation}
\mathsf{S}_{n}(z)\ :=\ \frac{f_{n}(z)-q}{\gamma^{n}}\;\underset{n\uparrow
\infty}{\longrightarrow}\;\text{some }\mathsf{S}(z)\ =:\ \sum_{j=0}^{\infty
}\nu_{j}z^{j},\qquad z\in D. \label{SchroederLimit}%
\end{equation}
Moreover, the convergence is uniform on each compact subsets of the interior
$D^{\circ}$ of $\,D.$\thinspace\ Furthermore, the function $\mathsf{S}$
restricted to the reals is the unique solution of the so-called
\emph{Schr\"{o}der functional equation}\/ (see, e.g., Kuczma (1968)
\cite[Theorem~6.1, p.137]{Kuczma1968}),
\begin{equation}
\mathsf{S}\!\left(  _{\!_{\!_{\,}}}f(s)\right)  =\,\gamma\,\mathsf{S}%
(s),\qquad0\leq s\leq1,
\end{equation}
satisfying
\begin{equation}
\mathsf{S}(q)=0\quad\text{and}\quad\lim_{s\rightarrow q}\mathsf{S}^{\prime
}(s)=1.
\end{equation}
As a consequence of (\ref{SchroederLimit}),
\begin{equation}
\lim_{n\uparrow\infty}\gamma^{-n}\,\mathbf{P}(Z_{n}=k)\ =\ \nu_{k\,},\qquad
k\geq1. \label{LocalLimit}%
\end{equation}
Consequently, in the Schr\"{o}der case, these extreme ($k$ is fixed) lower
deviation probabilities $\,\mathbf{P}(Z_{n}=k)$\thinspace\ are positive and
decay to $0$ with order $\,\gamma^{n}.$\thinspace\ On the other hand, the
characteristics $\,\alpha\in(0,\infty)$\thinspace\ describes the behavior of
the limiting quantities $w(x)$ and $\mathbf{P}(W\leq x)$ as $x\downarrow0$. In
fact, according to Biggins and Bingham (1993) \cite{BigginsBingham1993}, there
is a continuous, positive multiplicatively periodic function $V$ such that
\begin{equation}
x^{1-\alpha}\,w(x)\ =\ V(x)+o(1)\quad\text{as}\,\ x\downarrow0.
\label{DensityAsymp}%
\end{equation}
Dubuc (1971) \cite{Dubuc1971} has shown that the function $\,V\,\ $can be
replaced by a constant $\,V_{0}>0$\thinspace\ if and only if $\,$%
\begin{equation}
\mathsf{S}\!\left(  _{\!_{\!_{\,}}}\varphi(h)\right)  \ =\ K_{0\,}h^{-\alpha
},\qquad h\geq0, \label{embeddable}%
\end{equation}
for some constant $\,K_{0}>0,$\thinspace\ where $\,\varphi=\varphi_{W}%
$\thinspace\ denotes the Laplace function of $\,W,$%
\begin{equation}
\varphi_{W}(h)\ :=\ \mathbf{E}e^{-hW},\qquad h\geq0. \label{not.phi}%
\end{equation}
We mention that condition (\ref{embeddable}) is certainly fulfilled if
$\,Z$\thinspace\ is embeddable (see \cite[p.96]{AsmussenHering1983}) into a
continuous-time Galton-Watson process (as in the case of a geometric offspring
law, see Example~\ref{Ex.geom} below).\smallskip

Now we turn to the \emph{B\"{o}ttcher case.} Here $\,\mu\geq2$\thinspace
\ (recall Definition~\ref{D.type}). Clearly, opposed to (\ref{LocalLimit}),
extreme lower deviation probabilities disappear, even $\,\mathbf{P}(Z_{n}%
<\mu^{n})=0$\thinspace\ for all $n\geq1$. Evidently,
\begin{equation}
\mathbf{P}(Z_{n}=\mu^{n})\ =\ \mathbf{P}(Z_{n-1}=\mu^{n-1})\,p_{\mu}%
^{(\mu^{n-1})}. \label{evidently}%
\end{equation}
Hence,
\begin{equation}
\mathbf{P}(Z_{n}=\mu^{n})\ =\ \prod_{j=0}^{n-1}p_{\mu}^{(\mu^{j})}%
\ =\ \exp\Bigl[\frac{\mu^{n}-1}{\mu-1}\log p_{\mu}\Bigr].
\label{mu^n-probability}%
\end{equation}
Next, \thinspace$\mathbf{P}(Z_{n}=\mu^{n}+1)=\mathbf{P}(Z_{n-1}=\mu
^{n-1})\,\mu^{n-1}\,p_{\mu+1}\,p_{\mu}^{\mu^{n-1}-1}$. Thus, from
(\ref{evidently}),
\begin{equation}
\mathbf{P}(Z_{n}=\mu^{n}+1)\ =\ p_{\mu}^{-1}\,p_{\mu+1}\,\mu^{n-1}%
\,\mathbf{P}(Z_{n}=\mu^{n}).
\end{equation}
For simplification, consider for the moment the special case $\,p_{\mu+j}%
>0,$\thinspace\ $j\geq0.$\thinspace\ Then, as in the previous representation,
for fixed $\,k\geq0$\thinspace\ and some positive constants $C_{k\,},$%
\begin{equation}
\mathbf{P}(Z_{n}=\mu^{n}+k)\ \sim\ C_{k}\,\mu^{nk}\,\mathbf{P}(Z_{n}=\mu
^{n})\quad\text{as}\quad n\uparrow\infty. \label{general}%
\end{equation}
Consequently, in contrast to (\ref{LocalLimit}) in the Schr\"{o}der case, here
the lower positive deviation probabilities $\,\mathbf{P}(Z_{n}=\mu^{n}%
+k)$\thinspace\ do \emph{not}\/ have a uniform order of decay. But by
(\ref{general}),%
\begin{equation}
\mu^{-n}\log\mathbf{P}(Z_{n}=\mu^{n}+k)\;\underset{n\uparrow\infty
}{\longrightarrow}\;\log p_{\mu},\qquad k\geq0. \label{logarithm}%
\end{equation}
That is, on a \emph{logarithmic}\/ scale, we gain again a uniform order,
namely $\,-\mu^{n}.$

Turning back to the general B\"{o}ttcher case,%
\begin{equation}
\lim_{n\uparrow\infty}\left(  _{\!_{\!_{\,}}}f_{n}(s)\right)  ^{\!(\mu^{-n}%
)}\ =:\ \mathsf{B}(s),\qquad0\leq s\leq1,\label{not.B}%
\end{equation}
exists, is continuous, positive, and satisfies the \emph{B\"{o}ttcher
functional equation}\/ $\,$%
\begin{equation}
\mathsf{B}\!\left(  _{\!_{\!_{\,}}}f(s)\right)  =\,\mathsf{B}^{\mu}%
(s),\qquad0\leq s\leq1,
\end{equation}
with boundary conditions
\begin{equation}
\mathsf{B}(0)=0\quad\text{and}\quad\mathsf{B}(1)=1
\end{equation}
(see, e.g., Kuczma (1968) \cite[Theorem~6.9, p.145]{Kuczma1968}).

Recalling that $\,\mu\geq2,$\thinspace\ define $\beta\in(0,1)$ by
\begin{equation}
\mu\,=\,m^{\beta}. \label{not.beta}%
\end{equation}
According to \cite[Theorem~3]{BigginsBingham1993}, there exists a positive and
multiplicatively periodic function $V^{\ast}$ such that
\begin{equation}
-\log\mathbf{P}(W\leq x)\ =\ x^{-\beta/(1-\beta)}\,V^{\ast}(x)+o(x^{-\beta
/(1-\beta)})\quad\text{as}\,\ x\downarrow0. \label{BoettcherLeftTail}%
\end{equation}
If additionally $\,\log\varphi_{W}(h)\sim-\kappa h^{\beta}\,\ $as
$\,h\uparrow\infty$\thinspace\ for some constant $\,\kappa>0,$\thinspace\ then
by Bingham (1988) \cite[formula~(4)]{Bingham1988},
\begin{equation}
-\log\mathbf{P}(W\leq x)\ \sim\ \beta^{-1}(1-\beta)(\kappa\beta)^{1/(1-\beta
)}\,x^{-\beta/(1-\beta)}\quad\text{as}\,\ x\downarrow0.
\end{equation}

\subsection{Lower deviation probabilities in the literature\label{SS.lit}}

What else is known in the literature on lower deviation probabilities of
$\,Z\,?$\thinspace\ In the \emph{Schr\"{o}der case}\/ $\,(\,0<\alpha
<\infty),\,\ $Athreya and Ney (1970) \cite{AthreyaNey1970} proved that in case
of mash $\,d=1$\thinspace\ and $\,\mathbf{E}Z_{1}^{2}<\infty,$\thinspace\ for
every $\,\varepsilon\in(0,\eta),$\thinspace\ where
\begin{equation}
\eta\,:=\,m^{\alpha/(3+\alpha)}\,>\,1, \label{not.eta}%
\end{equation}
there exists a positive constant $C_{\varepsilon}$ such that for all
$k\geq1,$
\begin{equation}
\left\vert _{\!_{\!_{\,_{{}}}}}m^{n}\mathbf{P}(Z_{n}=k)-w(k/m^{n})\right\vert
\ \leq\ C_{\varepsilon}\,\frac{\eta^{-n}}{km^{-n}}+(\eta-\varepsilon)^{-n}.
\label{AN}%
\end{equation}
The estimate (\ref{AN}) allows to get some information on lower deviation
probabilities. Indeed, in the general Schr\"{o}der case, from
(\ref{DensityAsymp}),%
\begin{equation}
w(x)\ \asymp\ x^{\alpha-1}\quad\text{as}\,\ x\downarrow0 \label{DubucBounds}%
\end{equation}
(meaning that there are positive constants $\,C_{1}$\thinspace\ and $\,C_{2}%
$\thinspace\ such that $\,C_{1}$\thinspace$x^{\alpha-1}\leq w(x)\leq C_{2}%
$\thinspace$x^{\alpha-1},$\thinspace\ $0<x\leq1).$\thinspace\ Together with
(\ref{AN}) this implies
\begin{equation}
\mathbf{P}(Z_{n}=k_{n})\ =\ m^{-n}\,w(k_{n}/m^{n})\bigg[1+O\Bigl(\frac
{m^{\alpha n}}{k_{n}^{\alpha}\eta^{n}}+\frac{m^{(\alpha-1)n}}{k_{n}^{\alpha
-1}\,(\eta-\varepsilon)^{n}}\Bigr)\bigg]\quad\text{as}\,\ n\uparrow\infty.
\end{equation}
We want to show that in important special cases the $O$--expression is
actually an $o(1).$ Recalling the definition (\ref{not.eta}) of $\eta,$ one
easily verifies that $\,m^{\alpha n}/k_{n}^{\alpha}\eta^{n}\rightarrow
0$\thinspace\ (as $n\uparrow\infty)$ if and only if $\,k_{n}/m^{n(2+\alpha
)/(3+\alpha)}\rightarrow\infty.$\thinspace\ Concerning the second $O$-term, if
additionally $\,\alpha\leq1,$ then $\,m^{(\alpha-1)n}/k_{n}^{\alpha-1}\leq
1$\thinspace\ provided that $\,k_{n}\leq m^{n}.$\thinspace\ Hence, here
$\,m^{(\alpha-1)n}/\!\left(  k_{n}^{\alpha-1}\,(\eta-\varepsilon)^{n}\right)
$\thinspace\ converges to zero if $\,\eta-\varepsilon>1$. On the other hand,
if $\,\alpha>1$ and $\,k_{n}/m^{n(2+\alpha)/(3+\alpha)}\rightarrow\infty
$\thinspace\ (which we needed for the first term), then $\,m^{(\alpha
-1)n}/\!\left(  k_{n}^{\alpha-1}\,(\eta-\varepsilon)^{n}\right)  \rightarrow
0$\thinspace\ provided that additionally $\,\varepsilon\leq m^{\alpha
/(3+\alpha)}-m^{(\alpha-1)/(3+\alpha)}.$\thinspace\ Altogether, under the
assumptions in \cite{AthreyaNey1970},
\begin{equation}
\mathbf{P}(Z_{n}=k_{n})\ =\ m^{-n}\,w(k_{n}/m^{n})\left(  _{\!_{\!_{\,}}%
}1+o(1)\right)  \quad\text{as}\,\ n\uparrow\infty\label{AN1}%
\end{equation}
provided that both $\,k_{n}\leq m^{n}$\thinspace\ and $\,k_{n}/m^{n(2+\alpha
)/(3+\alpha)}\rightarrow\infty$.

In \cite{AthreyaNey1970} it is also mentioned that according to an unpublished
manuscript of S.~Karlin, in the Schr\"{o}der case, for each embeddable
processes $Z$ of finite second moment,%
\begin{equation}
\lim_{n\uparrow\infty}\,\frac{m^{\alpha n}}{k_{n}^{\alpha-1}}\ \mathbf{P}%
(Z_{n}=k_{n})\,\ \text{exists in }(0,\infty),\ \text{provided that}%
\,\ k_{n}=o(m^{n}). \label{unpubl}%
\end{equation}
In the present situation, as we remarked after (\ref{DensityAsymp}),
$\,w(x)\sim V_{0}\,x^{\alpha-1}$\thinspace\ as $\,x\downarrow0$\thinspace
\ with $\,V_{0}>0.$\thinspace\ Hence, from (\ref{unpubl}), for some constant
$\,C>0,$%
\begin{equation}
\mathbf{P}(Z_{n}=k_{n})\ \sim\ C\,m^{-n}\,w(k_{n}/m^{n})\quad\text{as}%
\,\ n\uparrow\infty,
\end{equation}
which is compatible with (\ref{AN1}).

Intuitively, the asymptotic behavior of lower deviation probabilities should
be more related to characteristics as $\,\alpha$\thinspace\ and $\,\beta
$\thinspace\ than to the tail of the offspring distribution. Thus one can
expect that it is possible to describe lower deviation probabilities
successfully without the second moment assumption used in
\cite{AthreyaNey1970}. Actually, in \cite[Theorem~1]{NeyVidyashankar2004} one
finds the following \emph{claim.}

Suppose $p_{0}=0$ and $\,\mathbf{E}Z_{1\!}\log Z_{1}<\infty$. Then there exist
positive constants $\,C_{1}<C_{2}$\thinspace\ such that for $k_{n}%
\rightarrow\infty$ with $k_{n}=O(m^{n})$ as $n\uparrow\infty,$%
\begin{equation}
C_{1}\ \leq\ \liminf_{n\uparrow\infty}\,\frac{\mathbf{P}(Z_{n}=k_{n})}{A_{n}%
}\ \leq\ \limsup_{n\uparrow\infty}\,\frac{\mathbf{P}(Z_{n}=k_{n})}{A_{n}%
}\ \leq\ C_{2}, \label{NV1}%
\end{equation}
where
\begin{equation}
A_{n}\ :=\ \left\{
\begin{array}
[c]{ll}%
p_{1}^{n}\,k_{n}^{\alpha-1} & \text{if}\,\ \alpha<1,\vspace{4pt}\\
\theta_{n}\,p_{1}^{n} & \text{if}\,\ \alpha=1,\vspace{4pt}\\
m^{-n} & \text{if}\,\ 1<\alpha\leq\infty,
\end{array}
\right.  \label{not.An}%
\end{equation}
and $\,\theta_{n}:=\left[  _{\!_{\!_{\,}}}n+1-\log k_{n}/\log m\right]
\!.$\thinspace\ Furthermore, if $k_{n}=m^{n-\ell_{n}}$ for natural numbers
$\,\ell_{n\,}=O(n)$\thinspace\ as $n\uparrow\infty,$ then
\begin{equation}
\lim_{n\uparrow\infty}\,A_{n}^{-1}\,\mathbf{P}(Z_{n}=k_{n})\ =:\ C_{\lim
}\ \,\text{exists in }\,(0,\infty). \label{NV2}%
\end{equation}

\subsection{Contradictions\label{SS.contr}}

Let us test that claim by an example which allows explicit calculations.

\begin{example}
[\textbf{Geometric offspring law}]\label{Ex.geom}Consider the offspring
generating function
\begin{equation}
f(s)\ =\ \frac{s}{m-(m-1)s}\ =\ \sum_{j=1}^{\infty}m^{-1}\,(1-m^{-1}%
)^{j-1}\,s^{j},\qquad0\leq s\leq1,
\end{equation}
(with mean $m>1)$. Obviously, here $\,q=0,$\thinspace\ $\gamma=m^{-1}%
,$\thinspace\ hence $\,\alpha=1.$\thinspace\ For the $n^{\mathrm{th}}$ iterate
one easily gets
\begin{equation}
f_{n}(s)\ =\ \frac{s}{m^{n}-(m^{n}-1)s}\ =\ \sum_{j=1}^{\infty}m^{-n}%
\,(1-m^{-n})^{j-1}\,s^{j}.
\end{equation}
Thus,
\begin{equation}
\mathbf{P}(Z_{n}=k)\ =\ m^{-n}\,(1-m^{-n})^{k-1}\ \leq\ m^{-n}, \label{thus}%
\end{equation}
for all $\,n,k\geq1$. On the other hand, since $\,p_{1}=m^{-1},$\thinspace\ by
claim (\ref{NV1}) there is a constant $C>0$ such that for the considered
$\,k_{n\,}$,%
\begin{equation}
\mathbf{P}(Z_{n}=k_{n})\ \geq\ C\,\theta_{n}\,m^{-n} \label{contr.1}%
\end{equation}
for $n$ large enough. If, for example, $k_{n}=m^{n/2}$ then $\theta
_{n}\rightarrow\infty,$ and (\ref{contr.1}) contradicts (\ref{thus}).
Consequently, the left-hand part of claim (\ref{NV1}) cannot be true in the
case \thinspace$\alpha=1.$\hfill$\Diamond$
\end{example}

Next we compare the claim with our discussion in the previous section on lower
deviation probabilities based on \cite{AthreyaNey1970}. In fact, under the
assumptions in \cite{AthreyaNey1970}, if additionally $\,k_{n}=o(m^{n}%
)$\thinspace\ but $\,\,k_{n}/m^{n(2+\alpha)/(3+\alpha)}$ $\rightarrow$
$\infty$\thinspace\ as $\,n\uparrow\infty,$\thinspace\ then by (\ref{AN1}) and
(\ref{DubucBounds}),%
\begin{equation}
\mathbf{P}(Z_{n}=k_{n})\ \asymp\ m^{-n}\,\Big(\frac{k_{n}}{m^{n}%
}\Big)^{\!\alpha-1}.
\end{equation}
Thus, in the case $\,1<\alpha<\infty$\thinspace\ we get $\,\mathbf{P}%
(Z_{n}=k_{n})=o(m^{-n})$\thinspace\ which contradicts the positivity of
$\,C_{\lim}$ in claim (\ref{NV2}), hence of $\,C_{1}$ in claim (\ref{NV1}).

Here is one more consideration. According to claim (\ref{NV1}), under
$\,1<\alpha\leq\infty,$%
\begin{equation}
\mathbf{P}(Z_{n}=k)\ \geq\ C\,m^{-n}%
\end{equation}
for all $\,k\in\lbrack m^{\varepsilon n},m^{(1-\varepsilon)n}],$%
\thinspace\ $\varepsilon\in(0,1/2),$\thinspace\ and all $n$ large enough. Here
and later, $\,C$\thinspace\ refers to a generic positive constant which might
change its value from place to place. Hence,
\begin{align}
\mathbf{E}Z_{n}^{-1}\  &  \geq\ \sum_{k=m^{\varepsilon n}}^{m^{(1-\varepsilon
)n}}k^{-1}\,\mathbf{P}(Z_{n}=k)\ \\
&  \geq\ C\,m^{-n}\sum_{k=m^{\varepsilon n}}^{m^{(1-\varepsilon)n}}%
k^{-1}\ =\ C\,(1-2\varepsilon)\,n\,m^{-n}\left(  _{\!_{\!_{\,}}}1+o(1)\right)
\quad\text{as}\quad n\uparrow\infty.\nonumber
\end{align}
But by Ney and Vidyashankar (2003) \cite[Theorem~1]{NeyVidyashankar2003},
$\,\mathbf{E}Z_{n}^{-1}$\thinspace\ is asymptotically equivalent to $m^{-n}$
(in the case $\,1<\alpha\leq\infty),$ getting one more
contradiction.\smallskip

Looking into details of the proof of \cite[Theorem~1]{NeyVidyashankar2004},
the following formulas are claimed to be true:%
\begin{align}
&  2\pi\,C_{\lim}\ =\ \label{C_lim}\\[4pt]
&  \left\{
\begin{array}
[c]{ll}%
{\displaystyle\sum\limits_{j\geq1}}
\nu_{j}\,w^{\ast j}(1),\vspace{6pt} & \alpha<1,\\%
{\displaystyle\int\nolimits_{\pi/m}^{\pi}}
\left[  _{\!_{\!_{\,_{{}}}}}\mathsf{S}\!\left(  _{\!_{\!_{\,}}}\psi(u)\right)
-\mathsf{S}\!\left(  _{\!_{\!_{\,}}}\psi(-u)\right)  \right]  \mathrm{d}u, &
\alpha=1,\\%
{\displaystyle\sum\limits_{\ell\geq0}}
m^{\ell}%
{\displaystyle\int\nolimits_{\pi/m}^{\pi}}
\left[  _{\!_{\!_{\,_{{}}}}}f_{\ell}\!\left(  _{\!_{\!_{\,}}}\psi(u)\right)
+f_{\ell}\!\left(  _{\!_{\!_{\,}}}\psi(-u)\right)  \right]  \mathrm{d}u+%
{\displaystyle\int\nolimits_{-\pi/m}^{\pi/m}}
\psi(u)\,\mathrm{d}u,\vspace{4pt} & 1<\alpha<\infty,\\%
{\displaystyle\int\nolimits_{-\pi/m}^{\pi/m}}
\psi(u)\,\mathrm{d}u, & \alpha=\infty,
\end{array}
\right. \nonumber
\end{align}
with $\,\mathsf{S}$\thinspace\ from (\ref{SchroederLimit}) and where
$\,\psi=\psi_{W}$\thinspace\ denotes the characteristic function of $\,W,$%
\begin{equation}
\psi_{W}(u)\ :=\ \mathbf{E}e^{iuW},\qquad u\in\mathbb{R}. \label{not.psi}%
\end{equation}
Recall that $C_{\lim}>0$\thinspace\ according to the claim. Now, if
$\,\alpha<1,$\thinspace\ the positiveness of $C_{\lim}$ is obvious from this
formula, since the density function $w$ is positive. But the point is that the
claim $\,C_{\lim}>0$\thinspace\ is \emph{not}\/ true in all other cases.

In fact, consider first the case $1<\alpha<\infty$. It is well-known that
$\psi$ solves the equation $\,$%
\begin{equation}
\psi(mu)\ =\ f\!\left(  _{\!_{\!_{\,}}}\psi(u)\right)  \!,\qquad
u\in\mathbb{R}, \label{equ.psi}%
\end{equation}
(e.g.\ \cite[formula~(6.1)]{AsmussenHering1983}). Iterating, we obtain $\,$%
\begin{equation}
\psi(m^{\ell}u)\ =\ f_{\ell}\!\left(  _{\!_{\!_{\,}}}\psi(u)\right)  ,\qquad
u\in\mathbb{R},\quad\ell\geq1. \label{iter.equ.psi}%
\end{equation}
Thus,
\begin{equation}
\int_{\pi/m}^{\pi}\left[  _{\!_{\!_{\,_{{}}}}}f_{\ell}\!\left(  _{\!_{\!_{\,}%
}}\psi(u)\right)  +f_{\ell}\!\left(  _{\!_{\!_{\,}}}\psi(-u)\right)  \right]
\mathrm{d}u\ =\ m^{-\ell}\int_{\pi m^{\ell-1}}^{\pi m^{\ell}}\left[
_{\!_{\!_{\,}}}\psi(u)+\psi(-u)\right]  \,\mathrm{d}u.
\end{equation}
Therefore,
\begin{align}
&  \bigg|\sum_{\ell\geq0}m^{\ell}\int_{\pi/m}^{\pi}\left[  _{\!_{\!_{\,_{{}}}%
}}f_{\ell}\!\left(  _{\!_{\!_{\,}}}\psi(u)\right)  +f_{\ell}\!\left(
_{\!_{\!_{\,}}}\psi(-u)\right)  \right]  \mathrm{d}u\bigg|\nonumber\\
&  \leq\ \int_{\pi/m}^{\infty}\left[  _{\!_{\!_{\,_{{}}}}}\left\vert
_{\!_{\!_{\,}}}\psi(u)\right\vert \!+\left\vert _{\!_{\!_{\,}}}\psi
(-u)\right\vert \right]  \mathrm{d}u, \label{39}%
\end{align}
which is finite, since in the Schr\"{o}der case (see, for example,
\cite{AthreyaNey1972}, p.83, Lemma~1),%
\begin{equation}
\left\vert _{\!_{\!_{\,}}}\psi(u)\right\vert \,\leq\ c\,|u|^{-\alpha},\qquad
u\in\mathbb{R}. \label{psi.est.Sch}%
\end{equation}
Hence,
\begin{equation}
\sum_{\ell\geq0}m^{\ell}\int_{\pi/m}^{\pi}\left[  _{\!_{\!_{\,_{{}}}}}f_{\ell
}\!\left(  _{\!_{\!_{\,}}}\psi(u)\right)  +f_{\ell}\!\left(  _{\!_{\!_{\,}}%
}\psi(-u)\right)  \right]  \mathrm{d}u\ =\ \biggl(\int_{-\infty}^{-\pi/m}%
+\int_{\pi/m}^{\infty}\biggr)\psi(u)\,\mathrm{d}u,
\end{equation}
and, consequently, by (\ref{C_lim}),
\begin{equation}
C_{\lim}\ =\ \frac{1}{2\pi}\int_{-\infty}^{\infty}\psi(u)\,\mathrm{d}u
\label{NV3}%
\end{equation}
in the present $\alpha\in(1,\infty)$ case. Inverting (\ref{not.psi}) gives
\begin{equation}
\int_{-\infty}^{\infty}e^{-iux}\,\psi(u)\,\mathrm{d}u\ =\ 2\pi w(x),\qquad
x>0.
\end{equation}
But by (\ref{DensityAsymp}) there is a (positive) constant $\,C$%
\thinspace\ such that $\,w(x)\leq C\,x^{\alpha-1},$\thinspace\ $0<x\leq
1.$\thinspace\ Hence, $\,w(0)=0,$\thinspace\ and (\ref{NV3}) implies $C_{\lim
}=0$.

In the case $\alpha=\infty,$ the proof of Lemma~5 in
\cite{NeyVidyashankar2004} is incorrect. In fact, the statement (82) there is
wrong. But we can start from (79) there\ (setting $\,\eta(r,s)\equiv1)\,\ $to
define%
\begin{equation}
I_{r-j}^{(2)}(r,s)\ :=\ \int_{\pi/m}^{\pi}e^{-ium^{-r+j}}f_{j}\!\left(
_{\!_{\!_{\,}}}\psi_{s+r-j}(u)\right)  \mathrm{d}u,\qquad r,s\geq1,\quad0\leq
j\leq r,
\end{equation}
where in this section by an abuse of notation,%
\begin{equation}
\psi_{\ell}(u)\,:=\,f_{\ell}(e^{iu/m^{\ell}})\,=\,\mathbf{E}e^{iuZ_{\ell
}/m^{\ell}},\qquad\ell\geq0,\quad u\in\mathbb{R}.
\end{equation}
By the global limit theorem (\ref{ILT}), for $\,u\in\mathbb{R}$\thinspace\ and
$\,j\geq0$\thinspace\ we get $\,\lim_{r,s\rightarrow\infty}\,\psi_{s+r-j}%
(u)$\ $=$\ $\psi(u)$\thinspace\ with $\,\psi=\psi_{W}$\thinspace\ from
(\ref{not.psi}), yielding $\,\lim_{r,s\rightarrow\infty}f_{j}\!\left(
_{\!_{\!_{\,}}}\psi_{s+r-j}(u)\right)  $\ $=$\ $f_{j}\!\left(  \psi
_{\!_{\!_{\,}}}(u)\right)  \!.$\thinspace\ Thus, by dominated convergence, for
$\,j\geq0,$
\begin{equation}
\lim_{r,s\rightarrow\infty}I_{r-j}^{(2)}(r,s)\ =\ \int_{\pi/m}^{\pi}%
f_{j}\!\left(  \psi_{\!_{\!_{\,}}}(u)\right)  \mathrm{d}u.
\end{equation}
Using this and the bound (81) there, one can easily verify that
\begin{equation}
\lim_{r,s\rightarrow\infty}\sum_{j=0}^{r}I_{r-j}^{(2)}(r,s)\ =\ \sum
_{j=0}^{\infty}m^{j}\int_{\pi/m}^{\pi}f_{j}\!\left(  \psi_{\!_{\!_{\,}}%
}(u)\right)  \mathrm{d}u.
\end{equation}
This gives for $\,C_{\lim}$\thinspace\ in the case $\,\alpha=\infty$%
\thinspace\ the same formula as written\ in (\ref{C_lim}) for the case
$\,1<\alpha<\infty.$\thinspace\ Now, instead of (\ref{psi.est.Sch}), in the
B\"{o}ttcher case we have%
\begin{equation}
\left\vert _{\!_{\!_{\,}}}\psi(u)\right\vert \,\leq\ e^{-Cu^{\beta}},\qquad
u\in\mathbb{R},
\end{equation}
for some constant $C,$ see \cite[Theorem~23]{Dubuc1971.AIF}. Therefore we get
again (\ref{39}) and (\ref{NV3}) also in the B\"{o}ttcher case. Finally, by
our Remark~\ref{R.w.at.0} below, $\,w(0)=0\,\ $and again we arrive at
$\,C_{\lim}=0.$

It remains to discuss the case $\alpha=1$. Here in the last formula at p.1156
of \cite{NeyVidyashankar2004} there is a sign error: It must be read as
$\,\int_{\pi/m}^{\pi}[\mathsf{S}\!\left(  _{\!_{\!_{\,}}}\psi(u)\right)
+\mathsf{S}\!\left(  _{\!_{\!_{\,}}}\psi(-u)\right)  ]\,\mathrm{d}%
u,$\thinspace\ which equals indeed the true value of $\,C_{\lim\,}.$%
\thinspace\ Now, at least if $Z$ is embeddable into a continuous-time
Galton-Watson process then analogously to (\ref{embeddable}) we get the
identity $\,\mathsf{S}\!\left(  _{\!_{\!_{\,}}}\psi(u)\right)  =K_{0}%
\,(iu)^{-1}\,\ $for some constant $\,K_{0}>0,$\thinspace\ implying
$\,\mathsf{S}\!\left(  _{\!_{\!_{\,}}}\psi(u)\right)  +\mathsf{S}\!\left(
_{\!_{\!_{\,}}}\psi(-u)\right)  \equiv0.$\thinspace\ Then $\,C_{\lim}%
=0$\thinspace\ for this class of processes.\smallskip

Altogether, all these contradictions to the quoted claim from
\cite[`Theorem~1']{NeyVidyashankar2004} (and its generalization
\cite[`Theorem~2']{NeyVidyashankar2004}) had been rather unexpected for us. Of
course, they gave us some more motivation to ask for the right and general
picture on lower deviation probabilities. Actually, it is wrong to distinguish
between velocity cases as in (\ref{not.An}). The only needed velocity case
differentiation is the mentioned dichotomy of Definition~\ref{D.dichot}. This
we will explain in the next two sections. In the end of Section~\ref{SS.Boe}
we then discuss the influence of \cite[`Theorem~1']{NeyVidyashankar2004} to
other results in \cite{NeyVidyashankar2004}.

\subsection{Lower deviations in the Schr\"{o}der case\label{SS.Schr}}

We start by stating our results on lower deviation probabilities in the
Schr\"{o}der case. Recall that here $\,\mu=0$\thinspace\ or $\,1.$

\begin{theorem}
[\textbf{Schr\"{o}der case}]\label{T.Schroeder}Let the offspring law be of the
Schr\"{o}der type and of type $(d,\mu)$. Then for all $\,k_{n}\equiv
\mu\ (\mathrm{mod}\,d)$ with $\,k_{n}\rightarrow\infty$\thinspace\ but
$\,k_{n}=o(c_{n}),$
\begin{equation}
\mathbf{P}(Z_{n}=k_{n})\ =\ \frac{d}{m^{n-a_{n}}\,c_{a_{n}}}\ w\Bigl(\frac
{k_{n}}{m^{n-a_{n}}\,c_{a_{n}}}\Bigr)\!\left(  _{\!_{\!_{\,}}}1+o(1)\right)
\label{SchroederAsymp}%
\end{equation}
and
\begin{equation}
\mathbf{P}(0<Z_{n}\leq k_{n})\ =\ \mathbf{P}\Bigl(0<W<\frac{k_{n}}{m^{n-a_{n}%
}\,c_{a_{n}}}\Bigr)\!\left(  _{\!_{\!_{\,}}}1+o(1)\right)
\label{SchroederAsymp1}%
\end{equation}
as\/$\,\ n\uparrow\infty,$\thinspace\ where for $\,n\geq1$\thinspace\ fixed we
put $\,a_{n}:=\min\{\ell\geq1:\,c_{\ell}\geq k_{n}\}$.
\end{theorem}

The appearing of the $\,a_{n}$\thinspace\ in the theorem, depending on the
$\,c_{n}$\thinspace\ and $\,k_{n}$\thinspace\ looks a bit disturbing, so we
have to discuss it. First assume additionally that $\,\mathbf{E}Z_{1\!}\log
Z_{1}<\infty.$\thinspace\ Since here we set $\,c_{n}=m^{n}$,\thinspace\ from
(\ref{SchroederAsymp}) we obtain the $a_{n}$-free formula%
\begin{equation}
\mathbf{P}(Z_{n}=k_{n})\ =\ d\,m^{-n}\,w(k_{n}/m^{n})\left(  _{\!_{\!_{\,}}%
}1+o(1)\right)  \!.
\end{equation}
Also, comparing this with (\ref{Intuition}), we see that under this
$\,Z_{1\!}\log Z_{1}$--moment condition in the Schr\"{o}der case,
$\,m^{-n}\,w(k_{n}/m^{n})$\thinspace\ describes not only normal deviation
probabilities but also lower ones.

On the other hand, without this additional moment condition, recalling
property (\ref{cn2}), $\,c_{n}=m^{n}\,L(m^{n})\,\ $with $\,L$\thinspace
\ slowly varying at infinity. Hence, we have
\begin{equation}
\frac{1}{m^{n-a_{n}}\,c_{a_{n}}}\ =\ \frac{1}{c_{n}}\,\frac{L(m^{n}%
)}{L(m^{a_{n}})},\quad\text{thus}\quad\frac{k_{n}}{c_{a_{n}}m^{n-a_{n}}%
}\ =\ \frac{k_{n}}{c_{n}}\,\frac{L(m^{n})}{L(m^{a_{n}})}\,.
\end{equation}
Therefore, from (\ref{SchroederAsymp}),%
\begin{equation}
\frac{c_{n}\,\mathbf{P}(Z_{n}=k_{n})}{d\,w(k_{n}/c_{n})}\ =\ \frac{L(m^{n}%
)}{L(m^{a_{n}})}\,\frac{w\!\left(  _{\!_{\!_{\,}}}k_{n}L(m^{n})/c_{n}%
L(m^{a_{n}})\right)  }{w(k_{n}/c_{n})}\left(  _{\!_{\!_{\,}}}1+o(1)\right)
\!.
\end{equation}
Using now (\ref{DensityAsymp}), we find
\begin{equation}
\frac{c_{n}\,\mathbf{P}(Z_{n}=k_{n})}{d\,w(k_{n}/c_{n})}\ =\ \Bigl(\frac
{L(m^{n})}{L(m^{a_{n}})}\Bigr)^{\!\alpha}\ \frac{V\!\left(  _{\!_{\!_{\,}}%
}k_{n}L(m^{n})/c_{n}L(m^{a_{n}})\right)  }{V(k_{n}/c_{n})}\left(
_{\!_{\!_{\,}}}1+o(1)\right)  \!. \label{RatioDev}%
\end{equation}
Next we want to expel the disturbing $\,a_{n}$ from this formula.

It is well-known (Seneta (1976) \cite[p.23]{Seneta1976}) that the regularly
varying function $\,x\mapsto xL(x)$\thinspace\ asymptotically equals a
(strictly) increasing, continuous, regularly varying function $\,x\mapsto
R(x):=xL_{1}(x)\,\ $with slowly varying $\,L_{1\,}.$\thinspace\ Hence,
$\,L(x)\sim L_{1}(x)$\thinspace\ as $x\uparrow\infty$. Using now
\cite[Lemma~1.3]{Seneta1976}, we conclude that the inverse function $R^{\ast}$
of $R$ equals $\,x\mapsto xL^{\ast}(x)$, where $L^{\ast}$ is again a slowly
varying function.

Put $\,x_{n}:=R^{\ast}(k_{n}).$\thinspace\ Then $\,k_{n}=x_{n}L_{1}(x_{n}%
)$\thinspace\ by the definition of $\,R^{\ast}.$ Recalling that $x_{n}%
=k_{n}L^{\ast}(k_{n})$, we get the identity
\begin{equation}
L^{\ast}(k_{n})\,L_{1}(x_{n})\,=\,1,\qquad n\geq1. \label{identity}%
\end{equation}
For $n$ fixed, define $\,b_{n}:=\min\left\{  _{\!_{\!_{\,}}}\ell
\geq1:\,m^{\ell}L_{1}(m^{\ell})\geq k_{n}\right\}  \!.$\thinspace\ Combined
with $\,x_{n}L_{1}(x_{n})$\ $=$\ $k_{n}$\thinspace\ we get%
\begin{equation}
m^{b_{n}}\,L_{1}(m^{b_{n}})\ \geq\ x_{n}L_{1}(x_{n})\ >\ m^{b_{n}-1}%
\,L_{1}(m^{b_{n}-1}).
\end{equation}
But $\,x\mapsto xL_{1}(x)$\thinspace\ is increasing, and the previous chain of
inequalities immediately gives%
\begin{equation}
m^{b_{n}}\,\geq\,x_{n}\,>\,m^{b_{n}-1}. \label{*}%
\end{equation}
By (\ref{cn2}),
\begin{equation}
c_{b_{n}+1}\ =\ m^{b_{n}+1}\,L(m^{b_{n}+1})\ =\ m\ \frac{L(m^{b_{n}+1})}%
{L_{1}(m^{b_{n}})}\ m^{b_{n}}\,L_{1}(m^{b_{n}})\ \geq\ k_{n}%
\end{equation}
for all $\,n$\thinspace\ sufficiently large. Here, in the last step we used
$\,m>1,$\thinspace\ that the slowly varying functions $\,L$\thinspace\ and
$\,L_{1}$\thinspace\ are asymptotically equivalent, and the definition of
$\,b_{n\,}.$\thinspace\ Now $\,c_{b_{n}+1}\geq k_{n}$\thinspace\ implies $\,$%
\begin{equation}
b_{n}+1\geq a_{n\,}, \label{first}%
\end{equation}
by the definition of $\,a_{n\,}.$\thinspace\ On the other hand, $\,$%
\begin{equation}
m^{a_{n}+1}\,L_{1}(m^{a_{n}+1})\ =\ m\,\frac{L_{1}(m^{a_{n}+1})}{L(m^{a_{n}}%
)}\ c_{a_{n}}\ \geq\ k_{n}%
\end{equation}
for all $\,n$\thinspace\ sufficiently large. Here, in the last step we used
the definition of $\,a_{n\,}.$\thinspace\ This gives%
\begin{equation}
a_{n}+1\geq b_{n\,}, \label{second}%
\end{equation}
by the definition of $b_{n\,}.$\thinspace\ Entering with (\ref{second}) and
(\ref{first}) into (\ref{*}), we get%
\begin{equation}
m^{a_{n}+1}\,\geq\,x_{n}\,>\,m^{a_{n}-2}\quad\text{for all}%
\,\ n\,\ \text{sufficiently large.}%
\end{equation}
Therefore, recalling (\ref{identity}),
\begin{equation}
L(m^{a_{n}})\ \sim\ L(x_{n})\ \sim\ L_{1}(x_{n})\ \sim\ \frac{1}{L^{\ast
}(k_{n})}\quad\text{as}\,\ n\uparrow\infty.
\end{equation}
Entering this into (\ref{RatioDev}) gives%
\begin{equation}
\frac{c_{n}\,\mathbf{P}(Z_{n}=k_{n})}{d\,w(k_{n}/c_{n})}\ =\ \left[
_{\!_{\!_{\,}}}L(m^{n})\,L^{\ast}(k_{n})\right]  ^{\alpha}\ \frac{V\!\left(
_{\!_{\!_{\,}}}k_{n\,}L(m^{n})\,L^{\ast}(k_{n})/c_{n}\right)  }{V(k_{n}%
/c_{n})}\left(  _{\!_{\!_{\,}}}1+o(1)\right)  \!, \label{RatioDev2}%
\end{equation}
which contains $\,L^{\ast}$\thinspace\ instead of the $\,a_{n\,}.$

Note also that such reformulation of (\ref{SchroederAsymp}) reminds the
classical Cram\'{e}r theorem (see, for example, Petrov (1975) \cite[\S VIII.2]%
{Petrov1975}) on large deviations for sums of independent random variables.
There the ratio of a tail probability of a sum of independent variables and
the corresponding normal law expression is considered. The crucial role in
Cram\'{e}r's theorem is played by the so-called Cram\'{e}r series
$\,\lambda(s):=\sum_{k=0}^{\infty}\lambda_{k}s^{k},$\thinspace\ where the
coefficients $\lambda_{k}$ depend on the cumulants of the summands. For the
lower deviation probabilities of supercritical Galton-Watson processes we have
a more complex situation: It is not at all clear, how to find the input data
$L,L^{\ast},V$ [entering into (\ref{RatioDev2})] based only on the knowledge
of the offspring generating function $f$.

It was already noted after (\ref{DensityAsymp}) that if $Z$ is embeddable into
a continuous-time Galton-Watson process then $\,V(x)\equiv V_{0\,}.$%
\thinspace\ Consequently, for embeddable processes, (\ref{RatioDev2}) takes
the slightly simpler form
\begin{equation}
\frac{c_{n}\,\mathbf{P}(Z_{n}=k_{n})}{d\,w(k_{n}/c_{n})}\ =\ \left[
_{\!_{\!_{\,}}}L(m^{n})\,L^{\ast}(k_{n})\right]  ^{\alpha}\left(
_{\!_{\!_{\,}}}1+o(1)\right)  \!. \label{RatioDev3}%
\end{equation}

On the other hand, if $\,V$ is not constant, the influence of this function on
the asymptotic behavior of the ratio $c_{n}\,\mathbf{P}(Z_{n}=k_{n}%
)/w(k_{n}/c_{n})$ is relatively small. Indeed, from continuity and
multiplicatively periodicity of $V(x)$ we see that $0<V_{1}\leq V(x)\leq
V_{2}<\infty$,\thinspace\ $x>0,$\thinspace\ for some constants $\,V_{1}%
,V_{2\,}.$\thinspace\ Therefore, from (\ref{RatioDev2}),%
\begin{gather}
\frac{V_{1}}{V_{2}}\left[  _{\!_{\!_{\,}}}L(m^{n})\,L^{\ast}(k_{n})\right]
^{\alpha}\left(  _{\!_{\!_{\,}}}1+o(1)\right)  \ \leq\ \frac{c_{n}%
\,\mathbf{P}(Z_{n}=k_{n})}{d\,w(k_{n}/c_{n})}\ \\
\leq\ \frac{V_{2}}{V_{1}}\left[  _{\!_{\!_{\,}}}L(m^{n})\,L^{\ast}%
(k_{n})\right]  ^{\alpha}\left(  _{\!_{\!_{\,}}}1+o(1)\right)  \!.\nonumber
\end{gather}
Note also that for many offspring distributions the bounds $V_{1}$ and $V_{2}$
may be chosen close to each other. This "near-constancy" phenomenon was
studied by Dubuc (1982) \cite{Dubuc1982} and by Biggins and Bingham (1991,
1993) \cite{BigginsBingham1991,BigginsBingham1993}.

\subsection{Lower deviations in the B\"{o}ttcher case\label{SS.Boe}}

Recall that $\,\mu\geq2$\thinspace\ in the B\"{o}ttcher case.

\begin{theorem}
[\textbf{B\"{o}ttcher case}]\label{T.Boettcher}Let the offspring law be of the
B\"{o}ttcher type and of type $(d,\mu)$. Then there exist positive constants
$B_{1}$ and $B_{2}$ such that for all\/ $\,k_{n}\equiv\mu^{n}\,(\mathrm{mod}%
\,d)$\thinspace\ with $\,k_{n}\geq\mu^{n}$\thinspace\ but $\,k_{n}=o(c_{n}),$%
\begin{subequations}
\begin{align}
-B_{1}\  &  \leq\ \liminf_{n\uparrow\infty}\mu^{b_{n}-n}\log\!\left[
_{\!_{\!_{\,}}}c_{n}\,\mathbf{P}(Z_{n}=k_{n})\right] \label{LB}\\
\  &  \leq\ \limsup_{n\uparrow\infty}\mu^{b_{n}-n}\log\!\left[  _{\!_{\!_{\,}%
}}c_{n}\,\mathbf{P}(Z_{n}=k_{n})\right]  \ \leq\ -B_{2}, \label{UB}%
\end{align}
where $\,b_{n}:=\min\{\ell:\,c_{\ell\,}\mu^{n-\ell}\geq2k_{n}\}$. The
inequalities remain true if one replaces $\,c_{n}\,\mathbf{P}(Z_{n}=k_{n})$ by
$\,\mathbf{P}(Z_{n}\leq k_{n})$.
\end{subequations}
\end{theorem}

Let us add at this place the following remark.

\begin{remark}
[\textbf{Behavior of }$w$ \textbf{at }$0$]\label{R.w.at.0}In analogy with
(\ref{DubucBounds}), in the B\"{o}ttcher case one has%
\begin{equation}
\log w(x)\ \asymp\ -x^{-\beta/(1-\beta)}\quad\text{as}\,\ x\downarrow0
\label{w.at.0}%
\end{equation}
with $\beta$\thinspace\ from (\ref{not.beta}). This can be shown using
techniques from the proof of Theorem~\ref{T.Boettcher}; see
Remark~\ref{R.w.at.0.details} below.\hfill$\Diamond$
\end{remark}

Our results in the B\"{o}ttcher case are much weaker than the results in the
Schr\"{o}der case: We got only logarithmic bounds. But this is not unexpected,
recall our discussion around (\ref{logarithm}).

Repeating arguments as we used to obtain (\ref{RatioDev2}), from
Theorem~\ref{T.Boettcher} we get
\begin{equation}
\frac{\log\!\left[  _{\!_{\!_{\,}}}c_{n}\,\mathbf{P}(Z_{n}=k_{n})\right]
}{(k_{n}/c_{n})^{-\beta/(1-\beta)}}\ \asymp\ -\Bigl[L^{\ast}(k_{n}/m^{\beta
n})\,L^{1/(1-\beta)}(m^{n})\Bigr]^{\beta}\quad\text{as\ }\,n\uparrow\infty,
\label{LB1}%
\end{equation}
where $L^{\ast}$ is such that $R_{1}(x):=x^{(1-\beta)}L(x)$ and $R_{2}%
(x):=x^{1/(1-\beta)}L^{\ast}(x)$ are asymptotic inverses, i.e. $R_{1}%
(R_{2}(x))\sim x$\thinspace\ and $R_{2}(R_{1}(x))\sim x$\thinspace\ as
$\,x\uparrow\infty$.

Taking into account (\ref{w.at.0}), we conclude that
\begin{equation}
\frac{\log\!\left[  _{\!_{\!_{\,}}}c_{n}\,\mathbf{P}(Z_{n}=k_{n})\right]
}{\log w(k_{n}/c_{n})}\ \asymp\ \Bigl[L^{\ast}(k_{n}/m^{\beta n}%
)\,L^{1/(1-\beta)}(m^{n})\Bigr]^{\beta}\quad\text{as\ }\,n\uparrow\infty.
\end{equation}
\smallskip

Let us continue our discussion of the paper \cite{NeyVidyashankar2004}. The
main reason to study there lower deviation probabilities is the application to
large deviation probabilities for the ratio $Z_{n+1}/Z_{n\,},$\thinspace
\ stated as Theorems~3 and 4 there. Using our Theorem~\ref{T.Schroeder}
(instead of `Theorem~1' there) in the proof of \cite[Theorem~3]%
{NeyVidyashankar2004} concerning large deviation probabilities in the
Schr\"{o}der case, one can easily verify that one needs only to change the
quantity $B$ in \cite[Theorem~3]{NeyVidyashankar2004} to be $-\log p_{1}$ for
\emph{all}\/ $\alpha\in(0,\infty),$\thinspace\ in order to get the right
picture. On the other hand, \cite[Theorem~4]{NeyVidyashankar2004} concerning
large deviation probabilities in the B\"{o}ttcher case is true as it is
stated, since `Theorem~1' was used only to show that%
\begin{equation}
\lim_{n\uparrow\infty}\,\frac{1}{k_{n}}\log\!\left[  _{\!_{\!_{\,}}}%
m^{n}\,\mathbf{P}(Z_{n}=k_{n})\right]  \ =\ 0\quad\text{if}\quad\frac{\mu^{n}%
}{k_{n}}\;\underset{n\uparrow\infty}{\longrightarrow}\;0, \label{only}%
\end{equation}
see \cite[p.1163]{NeyVidyashankar2004}. Recalling that $c_{n}=m^{n}$ and
$L(x)\equiv L^{\ast}(x)\equiv1$ under $\,\mathbf{E}Z_{1\!}\log Z_{1}<\infty,$
using our (\ref{LB1}), one obtains%
\begin{equation}
\frac{1}{k_{n}}\log\!\left[  _{\!_{\!_{\,}}}m^{n}\,\mathbf{P}(Z_{n}%
=k_{n})\right]  \ \asymp\ -\Bigl(\frac{m^{\beta n}}{k_{n}}\Bigr)^{\!1/(1-\beta
)}\quad\text{as\ }\,n\uparrow\infty.
\end{equation}
But $m^{\beta}=\mu$ by definition (\ref{not.beta}) of $\,\beta,$%
\thinspace\ and (\ref{only}) follows indeed.

\section{Cram\'{e}r transforms applied to Galton-Watson processes}

Our way to prove Theorems~\ref{T.Schroeder} and \ref{T.Boettcher} is based on
the well-known Cram\'{e}r method (see, e.g., \cite[Chapter~8]{Petrov1975}),
which was developed to study large deviations for sums of independent random
variables. A key in this method is the so-called \emph{Cram\'{e}r transform}\/
defined as follows. A random variable $X(h)$ is called a Cram\'{e}r transform
(with parameter $h\in\mathbb{R}$) of the random real variable $X$ if
\begin{equation}
\mathbf{E}e^{itX(h)}\ =\ \frac{\mathbf{E}e^{(h+it)X}}{\mathbf{E}e^{hX}%
}\,,\qquad t\in\mathbb{R}.
\end{equation}
Of course, this transformation is well-defined if $\,\mathbf{E}e^{hX}<\infty.$

In what follows, we will \emph{always assume}\/ that our offspring law
additionally satisfies $\,p_{0}=0.$\thinspace\ This condition is not crucial
but allows a bit simplified exposition of auxiliary results formulated in
Lemma~\ref{LargeDeviations} below and of the proof of
Theorem~\ref{T.Schroeder} in Section~\ref{SS.Schroed} (see also
Remark~\ref{R.case0} below).

\subsection{Basic estimates}

Fix an offspring law of type $(d,\mu).$\thinspace\ Let $\,n\geq1.$%
\thinspace\ Since $Z_{n}>0,$ the Cram\'{e}r transforms $Z_{n}(-h/c_{n})$ exist
for all $h\geq0.$ Clearly, $\,\mathbf{E}e^{itZ_{n}(-h/c_{n})}=f_{n}%
(e^{-h/c_{n}+it})/f_{n}(e^{-h/c_{n}}).$\thinspace\ We want to derive upper
bounds of $\,f_{n}(e^{-h/c_{n}+it})$\thinspace\ on $\,\left\{  t\in
\mathbb{R}:\ c_{n}^{-1}\pi d^{-1}\leq|t|\leq\pi d^{-1}\right\}  .$%
\thinspace\ \vspace{1pt}For this purpose, it is convenient to decompose the
latter set into $\,\bigcup_{j=1}^{n}J_{j}$\thinspace\ where%
\begin{equation}
J_{j}\ :=\ \left\{  t:\ c_{j}^{-1}\pi d^{-1}\leq|t|\leq c_{j-1}^{-1}\pi
d^{-1}\right\}  \!,\qquad j\geq1.
\end{equation}
To prepare for this, we start with the following generalization of
\cite[Lemma~2]{DubucSeneta1976}.

\begin{lemma}
[\textbf{Preparation}]\label{g.f.bound}Fix $\,\varepsilon\in(0,1).$%
\thinspace\ There exists $\,\theta=\theta(\varepsilon)\in(0,1)$ such that
\[
\bigl|f_{\ell}(e^{-h/c_{\ell}+it/c_{\ell}})\bigr|\,\leq\,\theta,\quad\ell
\geq0,\ \,h\geq0,\ \,t\in J_{\varepsilon}:=\left\{  t:\ \varepsilon\pi
d^{-1}\leq|t|\leq\pi d^{-1}\right\}  \!.
\]

\end{lemma}

\begin{proof}
Put $\,g_{h,t}(x):=e^{-hx+itx},$\thinspace\ \thinspace$h,x\geq0,$%
\thinspace\ $t\in\mathbb{R}.$\thinspace\ Evidently,
\begin{align}
\left\vert _{\!_{\!_{\,}}}g_{h,t}(x)-g_{h,t}(y)\right\vert \  &  =\ \left\vert
e^{-hx}(e^{itx}-e^{ity})+e^{ity}(e^{-hx}-e^{-hy})\right\vert \\
\  &  \leq\ |e^{itx}-e^{ity}|+|e^{-hx}-e^{-hy}|\ \leq\ \left(  _{\!_{\!_{\,}}%
}h+|t|\right)  |x-y|.\nonumber
\end{align}
It means that for $\,H\geq1$\thinspace\ and $\,T\geq\pi d^{-1}$\thinspace
\ fixed, $\,\mathcal{G}:=\left\{  g_{h,t};\ 0\leq h\leq H,\ |t|\leq T\right\}
$\thinspace\ is a family of uniformly bounded and equi-continuous functions on
$\,\mathbb{R}_{+\,}.$\thinspace\ Therefore, by (\ref{not.cn}),
\begin{equation}
f_{\ell}(e^{-h/c_{\ell}+it/c_{\ell}})\,=\,\mathbf{E}g_{h,t}(Z_{\ell}/c_{\ell
})\,\rightarrow\,\mathbf{E}g_{h,t}(W)\quad\text{as }\,\ell\uparrow\infty,
\label{UniformConv}%
\end{equation}
uniformly on $\mathcal{G}$\thinspace\ (see, e.g., Feller (1971)
\cite[Corollary in Chapter VIII, \S 1, p.252]{Feller1971.vol.II.2nd}). Since
$W>0$ has an absolutely continuous distribution, and $\,t\in J_{\varepsilon}%
$\thinspace\ implies $\,|t|\leq T,$%
\begin{equation}
\sup_{0\leq h\leq H,\ t\in J_{\varepsilon}}\left\vert \mathbf{E}%
e^{-hW+itW}\right\vert <1. \label{UniformBound}%
\end{equation}
{F}rom (\ref{UniformConv}) and (\ref{UniformBound}) it follows that there
exist $\delta_{1}\in(0,1)$ and $\ell_{0}$ such that
\begin{equation}
\sup_{0\leq h\leq H,\ t\in J_{\varepsilon}}\bigl|f_{\ell}(e^{-h/c_{\ell
}+it/c_{\ell}})\bigr|\ \leq\ \delta_{1\,},\qquad\ell>\ell_{0\,}.
\label{h<1-n>n_0-Bound}%
\end{equation}
On the other hand, $\,\bigcup_{\ell=0}^{\ell_{0}}\left\{  e^{-h/c_{\ell
}+it/c_{\ell}};\ h\geq0,\ t\in J_{\varepsilon}\right\}  $ is a subset of a
compact subset $\,K$\thinspace\ of the unit disc $\,D,$\thinspace\ where
$\,K$\thinspace\ does not contain the $d^{\mathrm{th}}$ roots of unity. Thus
for some $\delta_{2}\in(0,1),$%
\begin{equation}
\sup_{0\leq h\leq H,\ t\in J_{\varepsilon}}\bigl|f_{\ell}(e^{-h/c_{\ell
}+it/c_{\ell}})\bigr|\ \leq\ \delta_{2\,},\qquad\ell\leq\ell_{0\,}.
\label{h<1-n<n_0-Bound}%
\end{equation}
In fact, from Definition~\ref{D.type},%
\begin{equation}
f_{\ell}(z)\ =\ \sum_{j=0}^{\infty}\mathbf{P}(Z_{\ell}=\mu^{\ell}%
+jd)\,z^{\mu^{\ell}+jd},\qquad\ell\geq0,\quad z\in D, \label{rep.f}%
\end{equation}
implying%
\begin{equation}
\left\vert _{\!_{\!_{\,}}}f_{\ell}(z)\right\vert \ \leq\ \Big|\sum
_{j=0}^{\infty}\mathbf{P}(Z_{\ell}=\mu^{\ell}+jd)\,z^{jd}\Big|.
\end{equation}
But the latter sum equals $1$ if and only if $\,z$\thinspace\ is a
$d^{\mathrm{th}}$ root of unity, that is, if it is of the form $\,e^{2\pi
i/d}.$\thinspace\ 

Combining (\ref{h<1-n>n_0-Bound}) and (\ref{h<1-n<n_0-Bound}) gives the claim
in the lemma under the addition that $\,h\leq H.$\thinspace\ Consider now any
$h>H$. In this case
\begin{equation}
\bigl|f_{\ell}(e^{-h/c_{\ell}+it/c_{\ell}})\bigr|\ \leq\ f_{\ell
}(e^{-1/c_{\ell}}).
\end{equation}
By (\ref{not.cn}) we have $\,$%
\begin{equation}
f_{\ell}(e^{-h/c_{\ell}})\,=\,\mathbf{E}e^{-hZ_{\ell}/c_{\ell}}\,\rightarrow
\,\mathbf{E}e^{-hW}\,\in\,(0,1]\,\ \text{as}\,\ \ell\uparrow\infty,
\label{Lapl.conv}%
\end{equation}
uniformly in $h$ from compact subsets of $\,\mathbb{R}_{+\,}.$\thinspace\ In
particular,%
\begin{equation}
\sup_{\ell\geq1}f_{\ell}(e^{-1/c_{\ell}})\,<\,1. \label{88}%
\end{equation}
This completes the proof.
\end{proof}

The following lemma generalizes \cite[Lemma~3]{DubucSeneta1976}.

\begin{lemma}
[\textbf{Estimates on }$J_{1},\ldots J_{n}$]\label{L.J.1n}There are constants
$\,A>0$\thinspace\ and $\,\,\theta\in(0,1)$\thinspace\ such that for
$\,h\geq0,\,\ t\in J_{j\,},\,\ $and $\,1\leq j\leq n,$
\begin{equation}
\big|f_{n}(e^{-h/c_{n}+it})\big|\ \leq\ \left\{
\begin{array}
[c]{ll}%
A\,p_{1}^{n-j+1} & \text{in the Schr\"{o}der case,}\vspace{4pt}\\
\theta^{(\mu^{n-j+1})} & \text{in all cases.}%
\end{array}
\right.  \label{J.1n}%
\end{equation}

\end{lemma}

\begin{proof}
By (\ref{cn1}), we have $\,\varepsilon:=\inf_{\ell\geq1}c_{\ell-1}/c_{\ell}%
\in(0,1).$\thinspace\ If $\,t\in J_{j\,},$\thinspace\ $j\geq1,$\thinspace
\ then evidently,%
\begin{equation}
\pi d^{-1}\ \geq\ c_{j-1}\,|t|\ \geq\ c_{j-1}\,c_{j}^{-1}\pi d^{-1}%
\ \geq\ \varepsilon\pi d^{-1},
\end{equation}
hence $\,c_{j-1}t\in J_{\varepsilon\,}.$\thinspace\ Thus, by
Lemma~\ref{g.f.bound},
\begin{equation}
U\ :=\ \bigcup_{j=1}^{\infty}\left\{  _{\!_{\!_{\,_{{}}}}}f_{j-1}%
(e^{-h+it});\ \,h\geq0,\ \,t\in J_{j}\right\}  \ \subseteq\ \theta
D\quad\text{with }\,0<\theta<1.
\end{equation}
{F}rom the representation (\ref{rep.f}), $f_{\ell}(z)\leq|z|^{(\mu^{\ell})}$
for all $\,\ell\geq0$\thinspace\ and $\,|z|\leq1.$\thinspace\ Hence, for all
$\,z\in U\subseteq\theta D$\thinspace\ we have the bound $\left\vert
_{\!_{\!_{\,}}}f_{\ell}(z)\right\vert \leq\theta^{(\mu^{\ell})}.$%
\thinspace\ Thus, for $\,h\geq0,\ \,t\in J_{j\,},\,\ $and$\ \,1\leq j\leq n,$%
\begin{equation}
\bigl|f_{n}(e^{-h/c_{n}+it})\bigr|\ \leq\ f_{n-j+1}\left(  \bigl|f_{j-1}%
(e^{-h/c_{n}+it})\bigr|\right)  \ \leq\ \theta^{(\mu^{n-j+1})},
\label{iterate}%
\end{equation}
which is the second claim in (\ref{J.1n}).

If additionally $\,p_{1}>0,$ then by (\ref{SchroederLimit}) (and our
assumption $\,p_{0}=0)$\thinspace\ we have that $\,p_{1}^{-\ell}f_{\ell}%
(z)$\thinspace\ converges as $\,\ell\uparrow\infty,$\thinspace\ uniformly on
each compact $\,K\subset D^{\circ}.$\thinspace\ Therefore, there exists a
constant $C=C(K)$ such that
\begin{equation}
\left\vert _{\!_{\!_{\,}}}f_{\ell}(z)\right\vert \,\leq\,C\,p_{1\,}^{\ell
},\qquad\ell\geq0,\quad z\in K. \label{on.compact}%
\end{equation}
Consequently, iterating as in (\ref{iterate}),
\begin{equation}
\bigl|f_{n}(e^{-h/c_{n}+it})\bigr|\leq\ C\,p_{1}^{n-j+1},\qquad h\geq0,\quad
t\in J_{j\,},\quad1\leq j\leq n,
\end{equation}
finishing the proof.
\end{proof}

\subsection{On concentration functions}

Fix for the moment $\,h\geq0$\thinspace\ and $\,n\geq1.$\thinspace\ Denote by
$\left\{  _{\!_{\!_{\,}}}X_{j}(h,n)\right\}  _{j\geq1}$ \vspace{1pt}a sequence
of independent random variables which equal in law the Cram\'{e}r transform
$\,Z_{n}(-h/c_{n}),$\thinspace\ that is%
\begin{equation}
\mathbf{P}\!\left(  _{\!_{\!_{\,}}}X_{1}(h,n)=k\right)  \ =\ \frac
{e^{-kh/c_{n}}}{f_{n}(e^{-h/c_{n}})}\ \mathbf{P}(Z_{n}=k),\qquad k\geq1.
\label{CramerTransform}%
\end{equation}
Put $\,$%
\begin{equation}
S_{\ell}(h,n)\ :=\ \sum_{j=1}^{\ell}X_{j}(h,n),\qquad\ell\geq1.
\label{not.S_ell}%
\end{equation}
Note that%
\begin{equation}
\mathbf{E}e^{itS_{\ell}(h,n)}\ =\ \big(f_{n}(e^{-h/c_{n}+it})/f_{n}%
(e^{-h/c_{n}})\big)^{\ell}. \label{char.S_ell}%
\end{equation}
Recall notation $\alpha\in(0,\infty]$ from (\ref{not.gamma.alpha}).

\begin{lemma}
[\textbf{A concentration function estimate}]\label{ConFunction}For every
$\,h\geq0,$ there is a constant $\,A(h)$ such that
\begin{equation}
\sup_{n,k\geq1}c_{n}\,\mathbf{P}\!\left(  _{\!_{\!_{\,}}}S_{\ell
}(h,n)=k\right)  \ \leq\ \frac{A(h)}{\ell^{1/2}}\ ,\qquad\ell\geq\ell
_{0}:=1+[1/\alpha]. \label{not.l0}%
\end{equation}

\end{lemma}

\begin{proof}
It is known (see, for example, \cite[Lemma III.3, p.38]{Petrov1975}) that for
arbitrary (real-valued) random variables $X$ and every $\lambda,T>0,$
\begin{equation}
Q(X;\lambda)\ :=\ \sup_{y}\mathbf{P}(y\leq X\leq y+\lambda)\ \leq
\ \Bigl(\frac{96}{95}\Bigr)^{2}\max(\lambda,T^{-1})\int_{-T}^{T}\left\vert
_{\!_{\!_{\,}}}\psi_{X}(t)\right\vert \mathrm{d}t
\end{equation}
(with $\psi_{X}$ the characteristic function of $X)$. Applying this inequality
to $\,X=S_{\ell_{0}}(h,n)$\thinspace\ with $\,T=\pi d^{-1}$\thinspace\ and
$\,\lambda=1/2,$\thinspace\ using (\ref{char.S_ell}) we have%
\begin{equation}
\sup_{k\geq1}\mathbf{P}\!\left(  _{\!_{\!_{\,}}}S_{\ell_{0}}(h,n)=k\right)
\ \leq\ C\int_{-\pi d^{-1}}^{\pi d^{-1}}\frac{\left\vert f_{n}(e^{-h/c_{n}%
+it})\right\vert ^{\ell_{0}}}{f_{n}^{\ell_{0}}(e^{-h/c_{n}})}\,\mathrm{d}t
\end{equation}
for some constant $C$ independent of $\,h,n.$\thinspace\ By (\ref{Lapl.conv}),
for $\,h$\thinspace\ fixed, $f_{n}(e^{-h/c_{n}})$ is bounded away from zero,
and consequently, there is a positive constant $\,C(h)$\thinspace\ such that
\begin{equation}
\sup_{k\geq1}\mathbf{P}\!\left(  _{\!_{\!_{\,}}}S_{\ell_{0}}(h,n)=k\right)
\ \leq\ C(h)\int_{-\pi d^{-1}}^{\pi d^{-1}}\big|f_{n}(e^{-h/c_{n}%
+it})\big|^{\ell_{0}}\mathrm{d}t.\label{2.1}%
\end{equation}
Fist assume that $\alpha<\infty$ (Schr\"{o}der case). Using the first
inequality in (\ref{J.1n}), we get for $1\leq j\leq n,$
\begin{equation}
\int_{J_{j}}\big|f_{n}(e^{-h/c_{n}+it})\big|^{\ell_{0}}\,\mathrm{d}%
t\ \leq\ A^{\ell_{0}}\,p_{1}^{(n-j+1)\ell_{0}}|J_{j}|\ \leq\ 2\pi
d^{-1}A^{\ell_{0}}p_{1}^{(n-j+1)\ell_{0}}c_{j-1}^{-1}.\label{Jj}%
\end{equation}
On the other hand,
\begin{equation}
\int_{-\pi d^{-1}/c_{n}}^{\pi d^{-1}/c_{n}}\big|f_{n}(e^{-h/c_{n}%
+it})\big|^{\ell_{0}}\,\mathrm{d}t\ \leq\ 2\pi d^{-1}/c_{n\,}.\label{center}%
\end{equation}
{F}rom (\ref{Jj}) and (\ref{center}), for some constant $C,$
\begin{equation}
c_{n}\int_{-\pi d^{-1}}^{\pi d^{-1}}\big|f_{n}(e^{-h/c_{n}+it})\big|^{\ell
_{0}}\,\mathrm{d}t\ \leq\ C\Bigl(1+\sum_{j=1}^{n}p_{1}^{(n-j+1)\ell_{0}}%
c_{n}\,c_{j-1}^{-1}\Bigr).\label{both}%
\end{equation}
But by (\ref{cn1}), $\,$%
\begin{equation}
c_{n}\,\leq\,m^{n-j+1}c_{j-1\,},\qquad1\leq j\leq n.\label{cn1'}%
\end{equation}
Also, by the definition of $\ell_{0}$ in (\ref{not.l0}) and $\alpha$ in
(\ref{not.gamma.alpha}), $p_{1}^{\ell_{0}}m=p_{1}^{1+[1/\alpha]-1/\alpha}<1$.
Hence the right hand side of (\ref{both}) is bounded in $n$. Thus, from
(\ref{2.1}) it follows that
\begin{equation}
\sup_{n,k\geq1}c_{n}\,\mathbf{P}\!\left(  _{\!_{\!_{\,}}}S_{\ell_{0}%
}(h,n)=k\right)  \ \leq\ C(h).\label{l_0-bound}%
\end{equation}
This estimate actually holds also in the B\"{o}ttcher case, where $\,\ell
_{0}=1.$\thinspace\ Indeed, proceeding in the same way but using the second
inequality in (\ref{J.1n}) instead, the sum expression in (\ref{both}) has to
be replaced by%
\begin{equation}
\sum_{j=1}^{n}\theta^{(\mu^{n-j+1})}c_{n}\,c_{j-1}^{-1}\ \leq\ \sum_{j=1}%
^{n}\theta^{(\mu^{n-j+1})}\,m^{n-j+1}\ =\ \sum_{j=1}^{n}\theta^{(\mu^{j}%
)}\,m^{j},
\end{equation}
which again is bounded in $n.$

Note that (\ref{l_0-bound}) is (\ref{not.l0}) restricted to $\,\ell=\ell
_{0\,}.$\thinspace\ Hence, from now on we may restrict our attention to
$\,\ell>\ell_{0\,}.$\thinspace\ Let $Y_{1},\ldots,Y_{j}$ be independent
identically distributed random variables. Then by Kesten's inequality (see,
e.g., \cite[p.57]{Petrov1975}, there is a constant $\,C$\thinspace\ such that
for $\,0<\lambda^{\prime}<2\lambda$\thinspace\ the concentration function
inequality%
\begin{equation}
Q(Y_{1}+\ldots+Y_{j};\lambda)\ \leq\ \frac{C\lambda}{\lambda^{\prime}j^{1/2}%
}\ Q(Y_{1};\lambda)\bigl[1-Q(Y_{1};\lambda^{\prime})\bigr]^{-1/2}%
\end{equation}
holds. We specialize to $\,Y_{1}=S_{\ell_{0}}(h,n)$\thinspace\ and
$\,\lambda^{\prime}=\lambda=1/2.$\thinspace\ Note that $\,Q(Y_{1}%
;1/2)=\sup_{k\geq1}\mathbf{P}\!\left(  _{\!_{\!_{\,}}}S_{\ell_{0}%
}(h,n)=k\right)  <1$\thinspace\ in this case, since the random variable
$\,X_{1}(h,n)$\thinspace\ is non-degenerate. But also as $\,n\uparrow\infty
$\thinspace\ this quantity is bounded away from 1, which follows from
(\ref{l_0-bound}). Consequently, $\,\inf_{n\geq1}\!\left[  1-Q(Y_{1}%
;1/2)\right]  >0.$\thinspace\ Thus, using again (\ref{l_0-bound}), we infer
\begin{equation}
\sup_{n,k\geq1}\mathbf{P}\!\left(  _{\!_{\!_{\,}}}S_{j\ell_{0}}(h,n)=k\right)
\ \leq\ \frac{C_{1}(h)}{j^{1/2}}\ =\ \frac{C_{2}(h)}{(j\ell_{0})^{1/2}%
}\,,\qquad j\geq1, \label{jl_0-bound}%
\end{equation}
for some positive constants $C_{1}(h)$ and $C_{2}(h).$ If $X$ and $Y$ are
independent random variables, then, $Q(X+Y;\lambda)\leq Q(X;\lambda)$ (s.
\cite[Lemma III.1]{Petrov1975}). Thus for every $\ell>\ell_{0}$ we have the
inequality
\begin{equation}
\sup_{n,k\geq1}c_{n}\,\mathbf{P}\!\left(  _{\!_{\!_{\,}}}S_{\ell
}(h,n)=k\right)  \ \leq\ \sup_{n,k\geq1}c_{n}\,\mathbf{P}\!\left(
_{\!_{\!_{\,}}}S_{[\ell/\ell_{0}]\ell_{0}}(h,n)=k\right)  \!.
\end{equation}
Combining this bound once more with (\ref{jl_0-bound}), the proof is finished.
\end{proof}

\begin{remark}
[\textbf{Special case }$h=0$]\label{R.h0}Note that $\,S_{\ell}(0,n)$%
\thinspace\ equals in law to $\,Z_{n}$\thinspace\ conditioned to $\,Z_{0}%
=\ell.$\thinspace\ Therefore, by Lemma~\ref{ConFunction},%
\begin{equation}
\sup_{k\geq1}\mathbf{P}\!\left(  Z_{n}=k\,|\,Z_{0}=\ell\right)  \ \leq
\ \frac{A(0)}{\ell^{1/2}\,c_{n}}\,,\qquad n\geq1,\quad\ell\geq\ell_{0\,}.
\label{h.0}%
\end{equation}
In particular, if $\,\alpha>1,$\thinspace\ implying $\,\ell_{0}=1,$%
\thinspace\ in (\ref{h.0}) all initial states $\,Z_{0}$\thinspace\ are
possible. Especially, if $\,Z_{0}=1,$\thinspace\ then inequality (\ref{h.0})
generalizes the upper estimate in \cite[(10)]{NeyVidyashankar2004} to
processes without $Z_{1\!}\log Z_{1}$-moment condition.\hfill$\Diamond$
\end{remark}

Lemma~\ref{ConFunction} can also be used to get very useful bounds for
$\,\mathbf{P}\!\left(  Z_{n}=k\,|\,Z_{0}=\ell\right)  $\thinspace\ which are
not uniform in $k.$ This will be achieved in the next lemma by specializing
Lemma~\ref{ConFunction} to $\,h=1.$

\begin{lemma}
[\textbf{Non-uniform bounds}]\label{LargeDeviations}There exist two positive
constants $A$ and $\delta$ such that
\begin{equation}
c_{n}\,\mathbf{P}\!\left(  Z_{n}=k\,|\,Z_{0}=\ell\right)  \ \leq
\ A\,e^{k/c_{n}}\ell^{-1/2}\,e^{-\delta\ell},\qquad n,k\geq1,\quad\ell\geq
\ell_{0\,},
\end{equation}
\emph{[}with $\,\ell_{0}$\ defined in\/ \emph{(\ref{not.l0})].}
\end{lemma}

\begin{proof}
By the branching property and the definition (\ref{not.S_ell}) of $S_{\ell
}(h,n)$,%
\begin{equation}
\mathbf{P}\!\left(  Z_{n}=k\,|\,Z_{0}=\ell\right)  \ =\ e^{kh/c_{n}%
}\bigl[f_{n}(e^{-h/c_{n}})\bigr]^{\ell}\,\mathbf{P}\!\left(  _{\!_{\!_{\,}}%
}S_{\ell}(h,n)=k\right)  \!. \label{CramTrans}%
\end{equation}
Putting here $h=1$ and multiplying both sides by $\,c_{n\,},$\thinspace\ we
have
\begin{equation}
c_{n}\,\mathbf{P}\!\left(  Z_{n}=k\,|\,Z_{0}=\ell\right)  \ \leq\ e^{k/c_{n}%
}\bigl[f_{n}(e^{-1/c_{n}})\bigr]^{\ell}\max_{n,k\geq1}c_{n}\,\mathbf{P}%
\!\left(  _{\!_{\!_{\,}}}S_{\ell}(1,n)=k\right)  \!.
\end{equation}
Using Lemma~\ref{ConFunction} gives
\begin{equation}
c_{n}\,\mathbf{P}\!\left(  Z_{n}=k\,|\,Z_{0}=\ell\right)  \ \leq
\ A(1)\,\ell^{-1/2}\,e^{k/c_{n}}\bigl[f_{n}(e^{-1/c_{n}})\bigr]^{\ell}.
\label{3.1}%
\end{equation}
{F}rom (\ref{88}) the existence of a $\delta>0$ follows such that
$\,f_{n}(e^{-1/c_{n}})\leq e^{-\delta}$\thinspace\ for all $n\geq1$. Entering
this into (\ref{3.1}) finishes the proof.
\end{proof}

\subsection{On the limiting density function $w$}

Recall from Section~\ref{SS.growth} that $\,w$\thinspace\ denotes the density
function of $\,W,$\thinspace\ and $\,\psi=\psi_{W}$\thinspace\ its
characteristic function.

\begin{lemma}
[\textbf{Bounds for the limiting density}]\label{w(x)Bound} There is a
constant $\,A>0$\thinspace\ such that%
\begin{equation}
w^{\ast\ell}(x)\ \leq\ A\,\Bigl(\int_{0}^{x}w(t)\,\mathrm{d}t\Bigr)^{\ell
-\ell_{0}},\qquad x>0,\quad\ell\geq\ell_{0\,}.
\end{equation}

\end{lemma}

\begin{proof}
Suppose $\alpha<\infty$, the case $\alpha=\infty$ can be treated similarly. By
the inversion formula,
\begin{equation}
w^{\ast\ell_{0}}(x)\ =\ \frac{1}{2\pi}\int_{-\infty}^{\infty}e^{-itx}%
\,\psi^{\ell_{0}}(t)\,\mathrm{d}t,\qquad x>0.
\end{equation}
Hence,
\begin{equation}
A\ :=\ \sup_{x>0}w^{\ast\ell_{0}}(x)\ \leq\ \frac{1}{2\pi}\int_{-\infty
}^{\infty}\left\vert _{\!_{\!_{\,}}}\psi(t)\right\vert ^{\ell_{0}}\mathrm{d}t.
\label{l_0Density}%
\end{equation}
We want to convince ourselves that $\,A<\infty.$\thinspace\ For $\,j\geq0,$%
\begin{equation}
\int_{m^{j}}^{m^{j+1}}\left\vert _{\!_{\!_{\,}}}\psi(t)\right\vert ^{\ell_{0}%
}\mathrm{d}t\ =\ m^{j}\int_{1}^{m}\left\vert _{\!_{\!_{\,}}}\psi
(tm^{j})\right\vert ^{\ell_{0}}\mathrm{d}t\ =\ m^{j}\int_{1}^{m}\left\vert
f_{j}\!\left(  _{\!_{\!_{\,}}}\psi(t)\right)  \right\vert ^{\ell_{0}%
}\mathrm{d}t,
\end{equation}
where we used (\ref{iter.equ.psi}). Since $\,W>0$\thinspace\ has an absolute
continuous law, $\left\vert _{\!_{\!_{\,}}}\psi(t)\right\vert \leq C<1$ for
$t\in\lbrack1,m].$\thinspace\ Moreover, by (\ref{on.compact}), $\left\vert
_{\!_{\!_{\,}}}f_{j}(z)\right\vert \leq C\,p_{1}^{j}$ for $z$ in a compact
subset of $\,D^{\circ}.$\thinspace\ Therefore,
\begin{equation}
\int_{m^{j}}^{m^{j+1}}\left\vert _{\!_{\!_{\,}}}\psi(t)\right\vert ^{\ell_{0}%
}\mathrm{d}t\ \leq\ C\,m^{j}p_{1}^{j\ell_{0}}\ =\ C\,m^{j(1-\alpha\ell_{0})}%
\end{equation}
by definition (\ref{not.gamma.alpha}) of $\,\alpha.$\thinspace\ Consequently,
\begin{equation}
\int_{1}^{\infty}\left\vert _{\!_{\!_{\,}}}\psi(t)\right\vert ^{\ell_{0}%
}\mathrm{d}t\ \leq\ C\sum_{j=0}^{\infty}m^{j(1-\alpha\ell_{0})}\ <\ \infty,
\end{equation}
since $1-\alpha\ell_{0}<0$. Analogously,
\begin{equation}
\int_{-\infty}^{-1}\left\vert _{\!_{\!_{\,}}}\psi(t)\right\vert ^{\ell_{0}%
}\mathrm{d}t\ <\ \infty.
\end{equation}
Hence, $\,A$\thinspace\ in (\ref{l_0Density}) is finite. But $w^{\ast(\ell
+1)}(x)=\int_{0}^{x}w^{\ast\ell}(x-y)\,w(y)\,\mathrm{d}y,$\thinspace
\ $x>0,$\thinspace\ and the claim follows by induction.
\end{proof}

\subsection{A local central limit theorem}

Recall notation (\ref{not.S_ell}) of $\,S_{\ell}(h,n),$\thinspace\ $h\geq
0,$\thinspace\ $\ell,n\geq1.$\thinspace\ By an abuse of notation, denote by
$\psi_{\ell}=\psi_{\ell}^{h,n}$ the characteristic function of the random
variable
\begin{equation}
\ell^{-1/2}\,\sigma^{-1}(h,n)\left(  _{\!_{\!_{\,}}}S_{\ell}(h,n)-\mathbf{E}%
S_{\ell}(h,n)\right)  \!, \label{123}%
\end{equation}
where $\,\sigma(h,n):=\sqrt{\mathbf{E}\left(  _{\!_{\!_{\,}}}X_{1}%
(h,n)-\mathbf{E}X_{1}(h,n)\right)  ^{2}}.$\thinspace\ Note that by
(\ref{char.S_ell}),%
\begin{equation}
\psi_{\ell}^{h,n}(t)\ =\ \ \biggl(e^{-it\ell^{-1/2}\sigma^{-1}(h,n)\mathbf{E}%
X_{1}(h,n)}\,\frac{f_{n}(e^{-h/c_{n}+it\ell^{-1/2}\sigma^{-1}(h,n)})}%
{f_{n}(e^{-h/c_{n}})}\biggr)^{\!\ell}. \label{psi_ell}%
\end{equation}

\begin{lemma}
[\textbf{An Esseen type Inequality}]\label{EsseenInequality}If\/
$\,0<h_{1}\leq h_{2}<\infty$, then there exist positive constants\/
$C=C(h_{1},h_{2})$ and\/ $\varepsilon=\varepsilon(h_{1},h_{2})<1$%
\thinspace\ such that%
\begin{equation}
\sup_{h\in\lbrack h_{1},h_{2}],\ n\geq1}\,\bigl|\psi_{\ell}^{h,n}%
(t)-e^{-t^{2}/2}\bigr|\ \leq\ C\,\ell^{-1/2}\,|t|^{3}\,e^{-t^{2}/3}%
,\quad|t|<\varepsilon\,\ell^{1/2},\ \,\ell\geq1. \label{EssIn}%
\end{equation}

\end{lemma}

\begin{proof}
Put $\bar{X}_{j}(h,n):=X_{j}(h,n)-\mathbf{E}X_{j}(h,n)$. Using the global
limit theorem from (\ref{ILT}) one easily verifies that for some positive
constants $\,C_{1},\ldots,C_{4\,},$%
\begin{equation}
C_{1}\ \leq\ \frac{\sigma(h,n)}{c_{n}}\ \leq\ C_{2}\quad\text{uniformly in
}\,h\in\lbrack h_{1},h_{2}]\,\ \text{and }\,n\geq1 \label{VarBounds}%
\end{equation}
and%
\begin{equation}
C_{3}\ \leq\ \frac{\mathbf{E}\!\left\vert _{\!_{\!_{\,}}}\bar{X}%
_{1}(h,n)\right\vert ^{3}}{c_{n}^{3}}\ \leq\ C_{4}\quad\text{uniformly in
}\,h\in\lbrack h_{1},h_{2}]\,\ \text{and }\,n\geq1.
\end{equation}
Consequently, the Lyapunov ratio $\mathbf{E}\!\left\vert _{\!_{\!_{\,}}}%
\bar{X}_{1}(h,n)\right\vert ^{3}/\sigma^{3}(h,n)$ is bounded away from zero
and infinity. Applying now Lemma~V.1 from \cite{Petrov1975} to the random
variables $\,\bar{X}_{1}(h,n),\ldots,\bar{X}_{\ell}(h,n)$\thinspace\ we get
the desired result.
\end{proof}

The next lemma is a key step in our development concerning the B\"{o}ttcher
case. Recall notations $\,S_{\ell}:=S_{\ell}(h,n)$\thinspace\ and
$\,\sigma:=\sigma(h,n)$\thinspace\ defined in (\ref{not.S_ell}) and after
(\ref{123}), respectively.

\begin{lemma}
[\textbf{Local central limit theorem}]\label{LCLT}Suppose the offspring law is
of type $(d,\mu).$\thinspace\ If\/ $\,0<h_{1}\leq h_{2}<\infty$, then
\[
\sup_{\substack{h\in\lbrack h_{1},h_{2}]\\n\geq1}}\ \sup_{k:\,k\equiv\ell
\mu^{n}(\mathrm{mod}\,d)}\left\vert \ell^{1/2}\,\sigma(h,n)\,\mathbf{P}%
\!\left(  _{\!_{\!_{\,}}}S_{\ell}(h,n)=k\right)  -\frac{d}{\sqrt{2\pi}%
}\ e^{-x_{k,\ell}^{2}(h,n)/2}\right\vert \;\underset{\ell\uparrow\infty
}{\longrightarrow}\;0,
\]
where\/ $\,x_{k,\ell}:=x_{k,\ell}(h,n):=\ell^{-1/2}\,\sigma^{-1}(h,n)\left(
_{\!_{\!_{\,}}}k-\ell\,\mathbf{E}X_{1}(h,n)\right)  \!.$
\end{lemma}

Note that a local limit theorem, which would correspond to our case
$\,h=0$\thinspace\ but concerning an offspring law with finite variance and
with initial state tending to $\infty,$ was derived by H\"{o}pfner
\cite[Theorem~1]{Hoepfner1982}. The following proof of our lemma is a bit
simpler, since for $h>0$ the random variables $X_{1}(h,n)$ have finite moments
of all orders (also if the underlying $\,Z$\thinspace\ does not have finite
variance).\medskip

\noindent\emph{Proof of Lemma}\/ \ref{LCLT}.\thinspace$\ $By (\ref{char.S_ell}%
) and the inversion formula,%
\begin{equation}
\mathbf{P}\!\left(  _{\!_{\!_{\,}}}S_{\ell}=k\right)  \ =\ \frac{1}{2\pi}%
\int_{-\pi}^{\pi}e^{-itk}\Bigl[\frac{f_{n}(e^{-h/c_{n}+it})}{f_{n}%
(e^{-h/c_{n}})}\Bigr]^{\ell}\mathrm{d}t. \label{127}%
\end{equation}
Decomposing the unit circle,%
\begin{equation}
\left\{  e^{it}:\,-\pi<t\leq\pi\right\}  \ =\ \bigcup_{j=0}^{d-1}\left\{
\varrho^{j}\,e^{it}:\,-\pi d^{-1}<t\leq\pi d^{-1}\right\}  ,
\end{equation}
where $\,\varrho:=e^{2\pi i/d},$\thinspace\ the latter integral equals%
\begin{equation}
\sum_{j=0}^{d-1}\int_{-\pi d^{-1}}^{\pi d^{-1}}\varrho^{-jk}e^{-itk}%
\Bigl[\frac{f_{n}(\varrho^{j}e^{-h/c_{n}+it})}{f_{n}(e^{-h/c_{n}}%
)}\Bigr]^{\ell}\mathrm{d}t. \label{127'}%
\end{equation}
It is known (see, for instance, \cite[p.105]{AsmussenHering1983}) that for an
offspring law of type $(d,\mu)$ we have%
\begin{equation}
f_{n}(\varrho^{j}z)\ =\ \varrho^{j\mu^{n\!}}f_{n}(z),\qquad n,j\geq1,\quad
z\in D.
\end{equation}
Therefore the latter sum equals%
\begin{equation}
\int_{-\pi d^{-1}}^{\pi d^{-1}}e^{-itk}\Bigl[\frac{f_{n}(e^{-h/c_{n}+it}%
)}{f_{n}(e^{-h/c_{n}})}\Bigr]^{\ell}\mathrm{d}t\,\sum_{j=0}^{d-1}%
\varrho^{-j(k-\ell\mu^{n})}.
\end{equation}
But $\,\varrho^{-j(k-\ell\mu^{n})}\equiv1$\thinspace\ for $\,k\equiv\ell
\mu^{n}\ (\mathrm{mod}\ d).$\thinspace\ Altogether, for (\ref{127}) we get%
\begin{equation}
\mathbf{P}\!\left(  _{\!_{\!_{\,}}}S_{\ell}=k\right)  \ =\ \frac{d}{2\pi}%
\int_{-\pi d^{-1}}^{\pi d^{-1}}e^{-itk}\Bigl[\frac{f_{n}(e^{-h/c_{n}+it}%
)}{f_{n}(e^{-h/c_{n}})}\Bigr]^{\ell}\mathrm{d}t,\qquad k\equiv\ell\mu
^{n}\ (\mathrm{mod}\ d).
\end{equation}
Using the substitution $\,t\rightarrow t/\ell^{1/2}\sigma$\thinspace\ and
(\ref{psi_ell}), we arrive at%
\begin{equation}
\mathbf{P}\!\left(  _{\!_{\!_{\,}}}S_{\ell}=k\right)  \ =\ \frac{d}{2\pi
\ell^{1/2}\sigma}\int_{-\pi d^{-1}\ell^{1/2}\sigma}^{\pi d^{-1}\ell
^{1/2}\sigma}e^{-itx_{k,\ell}}\psi_{\ell}(t)\,\mathrm{d}t,\qquad k\equiv
\ell\mu^{n}\ (\mathrm{mod}\ d). \label{InvFormula}%
\end{equation}

Fix $\,0<h_{1}\leq h_{2}<\infty.$\thinspace\ Recall from (\ref{VarBounds})
that $\,$%
\begin{equation}
C_{1}\ \leq\ \inf_{h\in\lbrack h_{1},h_{2}],\,n\geq1}\,\frac{\sigma
(h,n)}{c_{n}}\ \leq\ \sup_{h\in\lbrack h_{1},h_{2}],\,n\geq1}\,\frac
{\sigma(h,n)}{c_{n}}\ \leq\ C_{2} \label{sup}%
\end{equation}
for some $\,0<C_{1}<C_{2}$\thinspace\ (depending on $\,h_{1},h_{2}%
).$\thinspace\ Choose a positive $\,$%
\begin{equation}
\varepsilon\,=\,\varepsilon(h_{1},h_{2})\,<\,C_{1}\pi d^{-1} \label{not.eps}%
\end{equation}
as in Lemma~\ref{EsseenInequality}. Take any $\,A=A(h_{1},h_{2})>\varepsilon
$\thinspace\ (to be specified later). Then the identity $\,\int_{-\infty
}^{\infty}e^{-itx-t^{2}/2}\,\mathrm{d}t=\sqrt{2\pi}\ e^{-x^{2}/2}$%
\thinspace\ and representation (\ref{InvFormula}) imply that
\begin{equation}
\sup_{k:\,k\equiv\ell\mu^{n\,}(\mathrm{mod}\,d)}\bigg|\ell^{1/2}%
\,\sigma\,\mathbf{P}\!\left(  _{\!_{\!_{\,}}}S_{\ell}=k\right)  -\frac
{d}{\sqrt{2\pi}}\ e^{-x_{k,\ell}^{2}/2}\bigg|\ \leq\ d\,(I_{1}+I_{2}%
+I_{3}+I_{4}),
\end{equation}
where
\begin{gather}
I_{1}\,:=\int_{-\varepsilon\ell^{1/2}}^{\varepsilon\ell^{1/2}}\bigl|\psi
_{\ell}(t)-e^{-t^{2}/2}\bigr|\,\mathrm{d}t,\quad I_{2}\,:=\int
_{|t|>\varepsilon\ell^{1/2}}e^{-t^{2}/2}\,\mathrm{d}t,\\
I_{3}\,:=\,\int_{\varepsilon\ell^{1/2}<|t|<A\ell^{1/2}}\left\vert
_{\!_{\!_{\,}}}\psi_{\ell}(t)\right\vert \mathrm{d}t,\quad I_{4}%
\,:=\int_{A\ell^{1/2}<|t|<\pi d^{-1}\ell^{1/2}\sigma}\!\left\vert
_{\!_{\!_{\,}}}\psi_{\ell}(t)\right\vert \mathrm{d}t.\nonumber
\end{gather}
[Of course, $I_{4}$ disappears if $\,A(h_{1},h_{2})>\pi d^{-1}\sigma(h,n).\,]$

Trivially, $I_{2}\rightarrow0$ as $\ell\uparrow\infty$. Further, due to
Lemma~\ref{EsseenInequality}, there is a $\,C=C(h_{1},h_{2})$\thinspace\ such
that
\begin{equation}
I_{1}\ \leq\ C\,\ell^{-1/2}\int_{0}^{\varepsilon\ell^{1/2}}t^{3}\,e^{-t^{2}%
/3}\,\mathrm{d}t\ \leq\ C\,\ell^{-1/2}\;\underset{\ell\uparrow\infty
}{\longrightarrow}\;0.
\end{equation}
Thus, it remains to show that the integrals $I_{3}$ and $I_{4}$ converge to
zero as $\ell\uparrow\infty,$\thinspace\ uniformly in the considered $h$ and
$n.$

First of all, using again (\ref{psi_ell}) and substituting $\,t\rightarrow
t\ell^{1/2}\sigma/c_{n\,},$\thinspace\ by (\ref{sup}) we obtain the following
estimates
\begin{subequations}
\label{I34}%
\begin{gather}
I_{3}\ \leq\ C_{2\,}\ell^{1/2}\int_{\varepsilon/C_{2\,}<|t|<A/C_{1}}%
\Big|\frac{f_{n}(e^{-h/c_{n}+it/c_{n}})}{f_{n}(e^{-h/c_{n}})}\Big|^{\ell
}dt,\label{I3}\\
I_{4}\ \leq\ C_{2\,}\ell^{1/2}\int_{A/C_{2}<|t|<\pi d^{-1}c_{n}}%
\Big|\frac{f_{n}(e^{-h/c_{n}+it/c_{n}})}{f_{n}(e^{-h/c_{n}})}\Big|^{\ell}dt.
\label{I4}%
\end{gather}

First we fix our attention to $I_{3\,}.$\thinspace\ By (\ref{UniformConv}),%
\end{subequations}
\begin{equation}
f_{n}(e^{-h/c_{n}+it/c_{n}})\,\rightarrow\,\mathbf{E}e^{-hW+itW}\quad\text{as
}\,n\uparrow\infty,
\end{equation}
uniformly in$\,\ h\in\lbrack0,h_{2}]$\thinspace\ and $\,t\in\lbrack0,A/C_{1}%
]$\thinspace\ [recall (\ref{not.eps})]. Also, by (\ref{Lapl.conv}),%
\begin{equation}
f_{n}(e^{-h/c_{n}})\,\rightarrow\,\mathbf{E}e^{-hW}\quad\text{as }%
\,n\uparrow\infty,\,\ \text{uniformly in}\,\ h\in\lbrack0,h_{2}]. \label{unif}%
\end{equation}
It follows that
\begin{equation}
\frac{f_{n}(e^{-h/c_{n}+it/c_{n}})}{f_{n}(e^{-h/c_{n}})}\;\underset
{n\uparrow\infty}{\longrightarrow}\;\frac{\mathbf{E}e^{-hW+itW}}%
{\mathbf{E}e^{-hW}}\ =\,\mathbf{E}e^{itW(-h)}, \label{RatioConv}%
\end{equation}
uniformly in$\,\ h\in\lbrack0,h_{2}]$\thinspace\ and $\,t\in\lbrack0,A/C_{1}%
]$\thinspace\ (with $\,W(-h)$\thinspace\ the Cram\'{e}r transform of
$\,W).$\thinspace\ Since the $\,W(-h)$\thinspace\ have absolutely continuous
laws, we have $\,|\mathbf{E}e^{itW(-h)}|<1$\thinspace\ for all $h\geq0$ and
$|t|>0$. This inequality and continuity of $\,(h,t)\mapsto\mathbf{E}%
e^{itW(-h)}$\thinspace\ imply that%
\begin{equation}
\sup_{0\leq h\leq h_{2\,},\ \varepsilon/C_{2}\leq|t|\leq A/C_{1}}%
\frac{|\mathbf{E}e^{-hW+itW}|}{\mathbf{E}e^{-hW}}\ <\ 1. \label{WBound}%
\end{equation}
Using (\ref{RatioConv}) and (\ref{WBound}) we infer the existence of a
positive constant $\eta=\eta(h_{1},h_{2})<1$\thinspace\ and an $\,n_{1}%
=n_{1}(h_{1},h_{2})\geq1$ such that for $n\geq n_{1\,},$%
\begin{equation}
\sup_{0\leq h\leq h_{2},\ \varepsilon/C_{2}\leq|t|\leq A/C_{1}}\Big|\frac
{f_{n}(e^{-h/c_{n}+it/c_{n}})}{f_{n}(e^{-h/c_{n}})}\Big|\ \leq\ \eta.
\label{previous}%
\end{equation}
Applying (\ref{previous}) to the bound of $\,I_{3}$\thinspace\ in (\ref{I3}),
we conclude that%
\begin{equation}
I_{3}\,\leq\,CA\,\ell^{1/2}\eta^{\ell}\,\rightarrow\,0\quad\text{as }%
\,\ell\uparrow\infty,
\end{equation}
uniformly in $h\in\lbrack h_{1},h_{2}]$ and $n\geq n_{1\,}.$\thinspace$\ $(The
remaining $n$ will be considered below.)

Next, we prepare for the estimation of $\,I_{4\,}.$\thinspace\ Since
$\,f_{n}(e^{-h/c_{n}})\geq f_{n}(e^{-h_{2}/c_{n}})$\thinspace\ for $\,0\leq
h\leq h_{2\,},$\thinspace\ and $\,f_{n}(e^{-h_{2}/c_{n}})\rightarrow
\mathbf{E}e^{-h_{2}W}>0$\thinspace\ as $\,n\uparrow\infty$\thinspace\ [recall
(\ref{unif})], there is a positive constant $C=C(h_{2})$ such that
\begin{equation}
\Big|\frac{f_{n}(e^{-h/c_{n}+it})}{f_{n}(e^{-h/c_{n}})}\Big|\ \leq
\ C\,\big|f_{n}(e^{-h/c_{n}+it})\big| \label{disting}%
\end{equation}
for all $t\in\mathbb{R},$\thinspace\ $0\leq h\leq h_{2\,},$\thinspace\ and
$\,n\geq1.$

At this point we have to distinguish between Schr\"{o}der and B\"{o}ttcher
cases. Actually, we proceed with the B\"{o}ttcher case $\alpha=\infty,$ which
is the only case we need later, and leave the other case for the reader.
Applying the second case of (\ref{J.1n}) to (\ref{disting}), we obtain the
estimate
\begin{equation}
\Big|\frac{f_{n}(e^{-h/c_{n}+it})}{f_{n}(e^{-h/c_{n}})}\Big|\ \leq
\ C\exp\!\left[  -\mu^{n-j+1}\log\theta^{-1}\right]  \!,
\end{equation}
$0\leq h\leq h_{2\,},\,\ t\in J_{j\,},\,\ $and $\,1\leq j\leq n.$%
\thinspace\ Since $\,\mu\geq2,$\thinspace\ there exists an $\,n_{2}%
=n_{2}(h_{2})$\thinspace\ such that
\begin{equation}
\Big|\frac{f_{n}(e^{-h/c_{n}+it})}{f_{n}(e^{-h/c_{n}})}\Big|\ \leq
\ \exp\!\left[  -\mu^{n-j}\log\theta^{-1}\right]  \!,
\end{equation}
if $\,0\leq h\leq h_{2\,},\,\ t\in J_{j\,},\,\ $and $\,1\leq j\leq n-n_{2\,}%
.$\thinspace\ But $\,|J_{j}|\leq2c_{j-1}^{-1}\pi d^{-1},$\thinspace\ hence%
\begin{equation}
\int_{J_{j}}\Big|\frac{f_{n}(e^{-h/c_{n}+it})}{f_{n}(e^{-h/c_{n}})}%
\Big|^{\ell}\,\mathrm{d}t\ \leq\ 2c_{j-1}^{-1}\pi d^{-1}\exp\!\left[
-\ell\,\mu^{n-j}\log\theta^{-1}\right]  \!.
\end{equation}
Summing over the considered $j$ gives%
\[
\int_{c_{n-n_{2}}^{-1}\pi d^{-1}\leq|t|\leq\pi d^{-1}}\Big|\frac
{f_{n}(e^{-h/c_{n}+it})}{f_{n}(e^{-h/c_{n}})}\Big|^{\ell}\,\mathrm{d}%
t\ \leq\ 2\pi d^{-1}\sum_{j=1}^{n-n_{2}}c_{j-1}^{-1}\exp\!\left[  -\ell
\,\mu^{n-j}\log\theta^{-1}\right]  \!,
\]
$0\leq h\leq h_{2}\,\ $and $\,n\geq n_{2\,}.$\thinspace\ Substituting
$t\rightarrow t/c_{n}$ and using (\ref{cn1'}), we arrive at
\begin{align}
&  \int_{\pi d^{-1}m^{n_{2}}\leq|t|\leq\pi d^{-1}c_{n}}\Big|\frac
{f_{n}(e^{-h/c_{n}+it/c_{n}})}{f_{n}(e^{-h/c_{n}})}\Big|^{\ell}dt\ \\
&  \leq\ 2\pi d^{-1}\sum_{j=1}^{n-n_{2}}m^{n-j+1}\exp\!\left[  -\ell
\,\mu^{n-j}\log\theta^{-1}\right]  \,\nonumber\\
&  \leq\ 2\pi d^{-1}\sum_{j=1}^{\infty}m^{j+1}\exp\!\left[  -\ell\,\mu^{j}%
\log\theta^{-1}\right]  \,\leq\ C\,e^{-C^{\prime}\ell}\nonumber
\end{align}
with constants \thinspace$C,C^{\prime},$\thinspace\ uniformly in
$\,h\in\lbrack h_{1},h_{2}]\,\ $and $\,n\geq n_{2\,}.$\thinspace\ Choosing now
$A$ so large that $\pi d^{-1}m^{n_{2}}\leq A/C_{2\,},$\thinspace\ we conclude
from (\ref{I4}) that
\begin{equation}
I_{4}\,\leq\,C\,\ell^{1/2}\,e^{-C^{\prime}\ell}\,\rightarrow\,0\quad\text{as
}\,\ell\uparrow\infty,
\end{equation}
uniformly in $\,h\in\lbrack h_{1},h_{2}]\,\ $and $\,n\geq n_{2\,}.$

Finally, we consider all $\,n\leq n^{\ast}:=n_{1}\vee n_{2\,}.$\thinspace\ By
definition, as in (\ref{rep.f}),%
\begin{equation}
\frac{f_{n}(e^{-h/c_{n}+it/c_{n}})}{f_{n}(e^{-h/c_{n}})}\ =\ \sum
_{j=0}^{\infty}\mathbf{P}\!\left(  _{\!_{\!_{\,}}}X_{1}(h,n)=\mu
^{n}+jd\right)  e^{(it/c_{n})(\mu^{n}+jd)}.
\end{equation}
Hence, since the set $\big\{e^{-it/c_{n}}:$\ $t\in\lbrack\varepsilon/C_{2},\pi
d^{-1}c_{n}]\big\}$ does not contain the $\,d^{\mathrm{th}}$\thinspace\ roots
of unity,
\begin{equation}
\sup_{t\in\lbrack\varepsilon/C_{2},\,\pi d^{-1}c_{n}]}\,\Big|\frac
{f_{n}(e^{-h/c_{n}+it/c_{n}})}{f_{n}(e^{-h/c_{n}})}\Big|\ =:\ \theta
_{n}(h)\ <\ 1. \label{1}%
\end{equation}
{F}rom the continuity $\,(h,t)\rightarrow f_{n}(e^{-h/c_{n}+it/c_{n}}%
)$\thinspace\ it follows that the function $\,\theta_{n}$\thinspace\ is
continuous, too. Therefore,%
\begin{equation}
\sup_{h\in\lbrack h_{1},h_{2}]}\,\theta_{n}(h)\,=:\,\bar{\theta}_{n}\,<\,1.
\label{2}%
\end{equation}
Combining (\ref{1}) and (\ref{2}),
\begin{equation}
\max_{n\leq n^{\ast}}\,\sup_{\substack{h\in\lbrack h_{1},h_{2}]\\t\in
\lbrack\varepsilon/C_{2},\,\pi d^{-1}c_{n}]}}\,\Big|\frac{f_{n}(e^{-h/c_{n}%
+it/c_{n}})}{f_{n}(e^{-h/c_{n}})}\Big|\ \leq\ \bar{\theta}%
\end{equation}
for some $\,\bar{\theta}<1.$\thinspace\ Substituting this into (\ref{I34})
gives
\begin{equation}
I_{3}+I_{4}\,\leq\,C\,\ell^{1/2}\,\bar{\theta}^{\ell}\,\rightarrow
\,0\quad\text{as }\,\ell\uparrow\infty,
\end{equation}
and the proof is finished.\hfill$\square$

\section{Proof of the main results}

\subsection{Schr\"{o}der case (proof of Theorem~\ref{T.Schroeder}%
)\label{SS.Schroed}}

Let $f,k_{n},$ and $a_{n}$ be as in Theorem~\ref{T.Schroeder}. Fix $\,n_{0}%
$\thinspace\ such that $\,c_{n}>k_{n}\geq1$\thinspace\ and $\,n>a_{n}\geq
1$\thinspace\ for all $\,n\geq n_{0\,},$\thinspace\ and consider only such
$\,n.$\thinspace\ Recall that $p_{0}=0$\thinspace\ by our convention. By the
Markov property,
\begin{equation}
\mathbf{P}(Z_{n}=k_{n})\ =\ \sum_{\ell=1}^{\infty}\mathbf{P}(Z_{n-a_{n}}%
=\ell)\,\mathbf{P}(Z_{a_{n}}=k_{n}\,|\,Z_{0}=\ell). \label{TotalProbability}%
\end{equation}
and
\begin{equation}
\mathbf{P}(Z_{n}\leq k_{n})\ =\ \sum_{\ell=1}^{\infty}\mathbf{P}(Z_{n-a_{n}%
}=\ell)\,\mathbf{P}(Z_{a_{n}}\leq k_{n}\,|\,Z_{0}=\ell).
\label{TotalProbability1}%
\end{equation}
\medskip

\noindent\emph{Step\ }$1^{\circ}$ (\emph{Proof of} (\ref{SchroederAsymp})).
\thinspace Using Lemma~\ref{LargeDeviations} we get for $N\geq\ell_{0}$ the
estimate
\begin{equation}
c_{a_{n}}\sum_{\ell=N}^{\infty}\mathbf{P}(Z_{n-a_{n}}=\ell)\,\mathbf{P}%
(Z_{a_{n}}=k_{n}\,|\,Z_{0}=\ell)\ \leq\ C\ \frac{e^{k_{n}/c_{a_{n}}}}{N^{1/2}%
}\ f_{n-a_{n}}(e^{-\delta}) \label{161}%
\end{equation}
for some constant $\,\delta>0.$\thinspace\ By (\ref{cn1}), and since
$\,c_{a_{n}-1}<k_{n}\leq c_{a_{n}}$\thinspace\ by the definition of
$\,a_{n\,},$%
\begin{equation}
m^{-1}\,\leq\,\frac{c_{a_{n}-1}}{c_{a_{n}}}\,\leq\,\frac{k_{n}}{c_{a_{n}}%
}\,\leq\,1. \label{4.1}%
\end{equation}
On the other hand, by (\ref{on.compact}),
\begin{equation}
f_{n-a_{n}}(e^{-\delta})\,\leq\,C\,p_{1}^{n-a_{n}}.
\end{equation}
Thus, from (\ref{161}),
\begin{equation}
p_{1}^{a_{n}-n}c_{a_{n}}\sum_{\ell=N}^{\infty}\mathbf{P}(Z_{n-a_{n}}%
=\ell)\,\mathbf{P}(Z_{a_{n}}=k_{n}\,|\,Z_{0}=\ell)\ \leq\ \frac{C}{N^{1/2}}\,.
\label{TailBound}%
\end{equation}

By \cite[Lemma~9]{DubucSeneta1976},
\begin{equation}
\lim_{n\uparrow\infty}\frac{1}{2\pi}\int_{-\pi d^{-1}c_{n}}^{\pi d^{-1}c_{n}%
}f_{n}^{\ell}(e^{it/c_{n}})\,e^{-itx}\,\mathrm{d}t\ =\ w^{\ast\ell}(x)
\end{equation}
uniformly in $x\in\lbrack m^{-1},1]$. This together with%
\begin{align}
&  c_{a_{n}}\mathbf{P}(Z_{a_{n}}=k_{n}\,|\,Z_{0}=\ell)\ \\
&  =\ \frac{d}{2\pi}\int_{-\pi d^{-1}c_{n}}^{\pi d^{-1}c_{n}}f_{a_{n}}^{\ell
}(e^{it/c_{n}})\,e^{-itk_{n}/c_{a_{n}}}\,\mathrm{d}t,\qquad\ell\equiv
k_{n}\ (\mathrm{mod}\ d),\nonumber
\end{align}
(see \cite[p.105]{AsmussenHering1983}) and (\ref{4.1}) gives
\begin{equation}
\lim_{n\uparrow\infty}\Bigl(c_{a_{n}}\mathbf{P}(Z_{a_{n}}=k_{n}\,|\,Z_{0}%
=\ell)-d\,w^{\ast\ell}(k_{n}/c_{a_{n}})\Bigr)\ =\ 0,\qquad\ell\equiv
k_{n}\ (\mathrm{mod}\ d). \label{Limit}%
\end{equation}
Since $k_{n}\equiv1\ (\mathrm{mod}\ d)$, the previous statement holds for all
$\ell\equiv1\ (\mathrm{mod}\ d)$. For other $\ell,$ the probabilities
$\mathbf{P}(Z_{n-a_{n}}=\ell)$ disappear. Thus, by (\ref{Limit}),
\begin{align}
&  \sum_{\ell=1}^{N-1}\mathbf{P}(Z_{n-a_{n}}=\ell)\,\mathbf{P}(Z_{a_{n}}%
=k_{n}\,|\,Z_{0}=\ell)\label{TruncatedSum}\\
\  &  =\ d\,c_{a_{n}}^{-1}\Bigl[\sum_{\ell=1}^{N-1}\mathbf{P}(Z_{n-a_{n}}%
=\ell)\,w^{\ast\ell}(k_{n}/c_{a_{n}})\Bigr]\left(  _{\!_{\!_{\,}}}%
1+o_{N}(1)\right) \nonumber
\end{align}
with $\,o_{N}(1)\rightarrow0\,\ $as $\,n\uparrow\infty,$\thinspace\ for each
fixed $\,N.$\thinspace\ Further, using Lemma~\ref{w(x)Bound}, one can easily
verify that there exist two constants $C$ and $\eta\in(0,1)$ such that
$\,w^{\ast\ell}(k_{n}/c_{a_{n}})\leq C\,\eta^{\ell}$\thinspace\ for all
$\,\ell\geq1$\thinspace\ and $\,n.$\thinspace\ Thus,%
\begin{equation}
\sum_{\ell=N}^{\infty}\mathbf{P}(Z_{n-a_{n}}=\ell)\,w^{\ast\ell}%
(k_{n}/c_{a_{n}})\ \leq\ C\sum_{\ell=N}^{\infty}\mathbf{P}(Z_{n-a_{n}}%
=\ell)\,\eta^{\ell}. \label{inequ}%
\end{equation}
But for every $\,\eta_{1}\in(\eta,1),$%
\begin{equation}
\sum_{\ell=N}^{\infty}\mathbf{P}(Z_{n-a_{n}}=\ell)\,\eta^{\ell}\ \leq
\ \Big(\frac{\eta}{\eta_{1}}\Big)^{\!N}f_{n-a_{n}}(\eta_{1})\ \leq
\ C\,\Big(\frac{\eta}{\eta_{1}}\Big)^{\!N}p_{1}^{n-a_{n}}, \label{N-bound}%
\end{equation}
where in the last step we used (\ref{on.compact}). Inequalities (\ref{inequ})
and (\ref{N-bound}) imply
\begin{equation}
\sum_{\ell=N}^{\infty}\mathbf{P}(Z_{n-a_{n}}=\ell)\,w^{\ast\ell}%
(k_{n}/c_{a_{n}})\ \leq\ C\,p_{1}^{n-a_{n}}\,e^{-\delta N} \label{TailBound2}%
\end{equation}
for all $\,n,N$\thinspace\ and some constant $\,\delta>0.$\thinspace
\ Combining (\ref{TotalProbability}), (\ref{TruncatedSum}), (\ref{TailBound})
and (\ref{TailBound2}), we have
\begin{gather}
\mathbf{P}(Z_{n}=k_{n})\ =\ d\,c_{a_{n}}^{-1}\Bigl[\sum_{\ell=1}^{\infty
}\mathbf{P}(Z_{n-a_{n}}=\ell)\,w^{\ast\ell}(k_{n}/c_{a_{n}})\Bigr]\left(
_{\!_{\!_{\,}}}1+o_{N}(1)\right) \label{4.2}\\
+\ O\!\left(  _{\!_{\!_{\,}}}c_{a_{n}}^{-1}\,p_{1}^{n-a_{n}}N^{-1/2}\right)
,\nonumber
\end{gather}
where the $O$-term applies to both $\,n,N\uparrow\infty.$\thinspace\ {By
}(\ref{iter.equ.psi}),%
\begin{equation}
m^{-j}\,w(x/m^{j})\ =\ \sum_{\ell=1}^{\infty}\mathbf{P}(Z_{j}=\ell
)\,w^{\ast\ell}(x),\qquad j\geq1,\quad x>0. \label{SelfSimilarity}%
\end{equation}
Putting here $\,j=n-a_{n\,},$ $x=k_{n}/c_{a_{n\,}},$\thinspace\ and
substituting into (\ref{4.2}), we arrive at%
\[
\mathbf{P}(Z_{n}=k_{n})\ =\ d\,c_{a_{n}}^{-1}\,m^{a_{n}-n}\,w(k_{n}m^{a_{n}%
-n}/c_{a_{n}})\left(  _{\!_{\!_{\,}}}1+o_{N}(1)\right)  +O\!\left(  c_{a_{n}%
}^{-1}\,p_{1}^{n-a_{n}}N^{-1/2}\right)  \!.
\]
{By} (\ref{DubucBounds}), (\ref{4.1}), and the definition
(\ref{not.gamma.alpha}) of $\,\alpha,$%
\begin{equation}
d\,c_{a_{n}}^{-1}\,m^{a_{n}-n}\,w(k_{n}m^{a_{n}-n}/c_{a_{n}})\ \geq
\ C\,c_{a_{n}}^{-1}\,m^{\alpha(a_{n}-n)}\ =\ C\,c_{a_{n}}^{-1}\,p_{1}%
^{n-a_{n}},\quad\text{for all }\,n.
\end{equation}
Therefore,
\begin{equation}
\mathbf{P}(Z_{n}=k_{n})\ =\ d\,c_{a_{n}}^{-1}\,m^{a_{n}-n}\,w(k_{n}m^{a_{n}%
-n}/c_{a_{n}})\left(  _{\!_{\!_{\,}}}1+o_{N}(1)+O(N^{-1/2})\right)  \!,
\end{equation}
where the $O$-term now applies to $\,N\uparrow\infty,$\thinspace\ uniformly in
$\,n.$\thinspace\ Letting first $\,n\uparrow\infty$\thinspace\ and then
$\,N\uparrow\infty,$\thinspace\ we see that (\ref{SchroederAsymp}) is
true.\medskip

\noindent\emph{Step\ }$2^{\circ}$ (\emph{Proof of} (\ref{SchroederAsymp1})).
\thinspace Trivially, for independent and identically distributed non-negative
random variables $\,X_{1},\ldots,\ X_{n}$\thinspace\ we have
\begin{equation}
\mathbf{P}(X_{1}+\ldots+X_{n}<x)\ \leq\ \mathbf{P}(\max_{j}X_{j}%
<x)\ =\ \mathbf{P}^{n}(X_{1}<x),\qquad x\geq0.
\end{equation}
Hence,
\begin{equation}
\mathbf{P}(Z_{a_{n}}\leq k_{n}\,|\,Z_{0}=\ell)\ \leq\ \mathbf{P}^{\ell
}(Z_{a_{n}}\leq k_{n}).
\end{equation}
Further, from (\ref{4.1}) and (\ref{ILT}),%
\begin{equation}
\mathbf{P}(Z_{a_{n}}\leq k_{n})\ \leq\ \mathbf{P}(c_{a_{n}}^{-1}Z_{a_{n}}%
\leq1)\;\underset{n\uparrow\infty}{\longrightarrow}\;\int_{0}^{1}w(x)\,dx.
\end{equation}
Therefore, since $w>0$ on all of $\,(0,\infty),$ there exists an $\eta
\in(0,1)$ such that $\,\mathbf{P}(Z_{a_{n}}\leq k_{n})\leq\eta$\thinspace\ for
all $\,n$\thinspace\ large enough. Thus,
\begin{equation}
\sum_{\ell=N}^{\infty}\mathbf{P}(Z_{n-a_{n}}=\ell)\,\mathbf{P}(Z_{a_{n}}\leq
k_{n}\,|\,Z_{0}=\ell)\ \leq\ \sum_{\ell=N}^{\infty}\mathbf{P}(Z_{n-a_{n}}%
=\ell)\,\eta^{\ell}%
\end{equation}
for all $\,N$\thinspace\ sufficiently large. Taking into account
(\ref{N-bound}), we conclude that
\begin{equation}
\sum_{\ell=N}^{\infty}\mathbf{P}(Z_{n-a_{n}}=\ell)\,\mathbf{P}(Z_{a_{n}}\leq
k_{n}\,|\,Z_{0}=\ell)\ \leq\ C\,p_{1}^{n-a_{n}}\,e^{-\delta N}
\label{TailBound3}%
\end{equation}
for $\,N$\thinspace\ sufficiently large and some $\,\delta>0.$\thinspace\ By
the same arguments,%
\begin{equation}
\sum_{\ell=N}^{\infty}\mathbf{P}(Z_{n-a_{n}}=\ell)\,F^{\ast\ell}%
(k_{n}/c_{a_{n}})\ \leq\ C\,p_{1}^{n-a_{n}}\,e^{-\delta N}, \label{TailBound4}%
\end{equation}
where $F(x):=\mathbf{P}(W<x),\,\ x\geq0.$

On the other hand, the continuity of $\,F$\thinspace\ and (\ref{ILT}) yield
that $\,\mathbf{P}(Z_{a_{n}}\leq c_{a_{n}}x\,|\,Z_{0}=\ell)\rightarrow
F^{\ast\ell}(x)$\thinspace\ uniformly in $\,x\geq0.$\thinspace\ Therefore,%
\begin{equation}
\lim_{n\uparrow\infty}\,\sup_{k\geq1}\Big|\mathbf{P}(Z_{a_{n}}\leq
k\,|\,Z_{0}=\ell)-F^{\ast\ell}(k/c_{a_{n}})\Big|\ =\ 0. \label{IntTheorem}%
\end{equation}

Combining (\ref{TotalProbability1}), (\ref{TailBound3}), (\ref{TailBound4}),
and (\ref{IntTheorem}), we arrive at
\begin{align}
&  \mathbf{P}(Z_{n}\leq k_{n})\ \label{178}\\
&  =\ \Bigl[\sum_{\ell=1}^{\infty}\mathbf{P}(Z_{n-a_{n}}=\ell)\,F^{\ast\ell
}(k_{n}/c_{a_{n}})\Bigr]\left(  _{\!_{\!_{\,}}}1+o_{N}(1)\right)
\,+\,O(p_{1}^{n-a_{n}}\,e^{-\delta N})\nonumber
\end{align}
with the same meaning of $\,o_{N}$\thinspace\ and the $O$-term as in the
previous step of proof. Since $\,\mathbf{P}(Z_{n-a_{n}}=1)=p_{1}^{n-a_{n}}%
$\thinspace\ and $\,F(k_{n}/c_{a_{n}})\geq F(m^{-1})>0$\thinspace\ by
(\ref{4.1}), we obtain%
\begin{equation}
p_{1}^{n-a_{n}}\,e^{-\delta N}\ \leq\ C\,e^{-\delta N}\sum_{\ell=1}^{\infty
}\mathbf{P}(Z_{n-a_{n}}=\ell)\,F^{\ast\ell}(k_{n}/c_{a_{n}}).
\end{equation}
Combining this inequality with (\ref{178}) gives
\begin{equation}
\mathbf{P}(Z_{n}\leq k_{n})\ =\ \Bigl[\sum_{\ell=1}^{\infty}\mathbf{P}%
(Z_{n-a_{n}}=\ell)\,F^{\ast\ell}(k_{n}/c_{a_{n}})\Bigr]\left(  _{\!_{\!_{\,}}%
}1+o_{N}(1)+O(e^{-\delta N})\right)  \!.
\end{equation}
Integrating both parts of (\ref{SelfSimilarity}), one has
\begin{equation}
F(y/m^{k})\ =\ \sum_{\ell=1}^{\infty}\mathbf{P}(Z_{k}=\ell)\,F^{\ast\ell
}(y),\qquad k\geq1,\quad y>0.
\end{equation}
Thus,%
\begin{equation}
\mathbf{P}(Z_{n}\leq k_{n})\ =\ F\Big(\frac{k_{n}}{c_{a_{n}}\,m^{n-a_{n}}%
}\Big)\left(  _{\!_{\!_{\,}}}1+o_{N}(1)+O(e^{-\delta N})\right)  \!.
\end{equation}
Letting again first $\,n\uparrow\infty$\thinspace\ and then $\,N\uparrow
\infty$\thinspace\ finishes the proof.\hfill$\square$

\begin{remark}
[\textbf{Proof in the case }$p_{0}>0$]\label{R.case0}We indicate now how to
proceed with the proof of Theorem~\ref{T.Schroeder} in the remaining case
$\,p_{0}>0.$\thinspace\ Here in the representation (\ref{TotalProbability})
one has additionally to take into account that%
\begin{align}
&  \mathbf{P}(Z_{a_{n}}=k_{n}\,|\,Z_{0}=\ell)\ \\
&  =\ \sum_{j=1}^{\ell}\left(
\genfrac{}{}{0pt}{1}{\ell}{j}%
\right)  f_{a_{n}}^{\ell-j}(0)\left(  _{\!_{\!_{\,}}}1-f_{a_{n}}(0)\right)
^{j}\,\mathbf{P}\!\Big\{\sum_{i=1}^{j}Z_{a_{n}}^{(i)}=k_{n}\;\Big|\;Z_{a_{n}%
}^{(i)}>0,\ 1\leq i\leq j\Big\},\nonumber
\end{align}
where the $\,Z^{(1)},Z^{(2)},\ldots$\thinspace\ are independent copies of
$\,Z.$\thinspace\ Then instead of Lemma~\ref{LargeDeviations} we need%
\[
c_{n}\,\mathbf{P}\!\Big\{\sum_{i=1}^{j}Z_{a_{n}}^{(i)}=k_{n}\;\Big|\;Z_{a_{n}%
}^{(i)}>0,\ 1\leq i\leq j\Big\}\,\leq\,A\,e^{k/c_{n}}j^{-1/2}\,e^{-\delta\ell
},\quad n,k\geq1,\ \,j\geq\ell_{0\,}.
\]
But this is valid by%
\begin{equation}
\mathbf{P}\big\{z^{Z_{n}^{(1)}}\,\big|\,Z_{n}^{(1)}>0\big\}\ =\ \frac
{f_{n}(z)-f_{n}(0)}{1-f_{n}(0)}\;\underset{n\uparrow\infty}{\longrightarrow
}\;\frac{\mathsf{S}(z)-\mathsf{S}(0)}{1-q}\,,
\end{equation}
uniformly in $\,z$\thinspace\ from compact subsets of $\,D^{\circ}.$%
\thinspace\ This indeed follows from (\ref{SchroederLimit}).\hfill$\Diamond$
\end{remark}

\subsection{B\"{o}ttcher case (proof of Theorem~\ref{T.Boettcher})}

From the Markov property,
\begin{equation}
\mathbf{P}(Z_{n}=k_{n})\ =\ \sum_{\ell=\mu^{n-b_{n}}}^{\infty}\mathbf{P}%
(Z_{n-b_{n}}=\ell)\,\mathbf{P}(Z_{b_{n}}=k\,|\,Z_{0}=\ell).
\label{BoettcherTP}%
\end{equation}
Using (\ref{CramTrans}) and Lemma~\ref{ConFunction}, we obtain the following
estimate
\begin{equation}
c_{b_{n}}\mathbf{P}(Z_{b_{n}}=k_{n}\,|\,Z_{0}=\ell)\ \leq\ A(h)\,\ell
^{-1/2}\bigl[e^{hk_{n}/\ell c_{b_{n}}}f_{b_{n}}(e^{-h/c_{b_{n}}})\bigr]^{\ell
}.
\end{equation}
{F}rom the definition of $\,b_{n}$ it immediately follows that
\begin{equation}
2k_{n}\ \leq\ c_{b_{n}}\mu^{n-b_{n}}\ =\ c_{b_{n}-1}\mu^{n-b_{n}+1}%
\Bigl(\frac{c_{b_{n}}}{\mu c_{b_{n}-1}}\Bigr)\ \leq\ 2k_{n}\frac{m}{\mu}\,.
\label{a_nbounds}%
\end{equation}
Hence,
\begin{equation}
\frac{hk_{n}}{\ell c_{b_{n}}}\ \leq\ \frac{h}{2} \label{198}%
\end{equation}
for \thinspace$\ell\geq\mu^{n-b_{n}}$. Therefore,
\begin{equation}
c_{b_{n}}\mathbf{P}(Z_{b_{n}}=k_{n}\,|\,Z_{0}=\ell)\ \leq\ A(h)\,\ell
^{-1/2}\bigl[e^{h/2}f_{b_{n}}(e^{-h/c_{b_{n}}})\bigr]^{\ell}.
\end{equation}
It is known (see, for example, \cite{AsmussenHering1983}, Corollary III.5.7),
that $\mathbf{E}W=1$ if $\mathbf{E}Z_{1\!}\log Z_{1}<\infty$ and
$\mathbf{E}W=\infty$ otherwise. It means, that for the Laplace function
$\,\varphi=\varphi_{W}$\thinspace\ of $\,W$\thinspace\ we have $\,e^{h/2}%
\varphi(h)<1$\thinspace\ for all small enough $h$. Thus, due to the global
limit theorem (\ref{ILT}), there exist $\delta<1$ and $h_{0}>0$ such that
$e^{h_{0}/2}f_{n}(e^{-h_{0}/c_{n}})\leq e^{-\delta}$ for all large enough $n$.
Hence,%
\begin{equation}
c_{b_{n}}\mathbf{P}(Z_{b_{n}}=k_{n}\,|\,Z_{0}=\ell)\ \leq\ A\,\ell
^{-1/2}\,e^{-\delta\ell}. \label{FinalBound}%
\end{equation}
Inserting (\ref{FinalBound}) into (\ref{BoettcherTP}), we obtain
\begin{equation}
c_{b_{n}}\mathbf{P}(Z_{n}=k_{n})\ \leq\ A\,\mu^{-(n-b_{n})/2}f_{n-b_{n}%
}(e^{-\delta}),
\end{equation}
consequently,
\begin{equation}
\mu^{b_{n}-n}\log\!\left[  _{\!_{\!_{\,}}}c_{n}\,\mathbf{P}(Z_{n}%
=k_{n})\right]  \ \leq\ \mu^{b_{n}-n}C+\mu^{b_{n}-n}\log\Bigl(\frac{c_{n}%
}{c_{b_{n}}}\Bigr)+\frac{\log f_{n}(e^{-\delta})}{\mu^{n-b_{n}}}\,.
\end{equation}
Since $c_{n}/c_{b_{n}}\leq m^{n-b_{n}}$ and $\mu^{n-b_{n}}=m^{\beta(n-b_{n})}%
$, $\mu^{b_{n}-n}\log(c_{n}/c_{b_{n}})\rightarrow0$ as $n\uparrow\infty$.
Thus,
\begin{equation}
\limsup_{n\uparrow\infty}\mu^{b_{n}-n}\log\!\left[  _{\!_{\!_{\,}}}%
c_{n}\,\mathbf{P}(Z_{n}=k_{n})\right]  \ \leq\ \limsup_{n\uparrow\infty}%
\frac{\log f_{n-b_{n}}(e^{-\delta})}{\mu^{n-b_{n}}}\,.
\end{equation}
Using (\ref{not.B}), we arrive at the desired upper bound.

We show now that (\ref{UB}) holds for $\,\log\mathbf{P}(Z_{n}\leq k_{n}%
).$\thinspace\ First of all we note that for arbitrary non-negative random
variable $X$ and all $x,h\geq0$
\begin{equation}
\mathbf{P}(X\leq x)\,\leq\,e^{hx}\,\mathbf{E}e^{-hX}.
\end{equation}
Applying this bound to the process $Z$ starting from $\ell$ individuals and
taking into account (\ref{198}), we have
\begin{equation}
\mathbf{P}(Z_{b_{n}}\leq k_{n}\,|\,Z_{0}=\ell)\ \leq\ \bigl[e^{hk_{n}/\ell
c_{b_{n}}}f_{b_{n}}(e^{-h/c_{b_{n}}})\bigr]^{\ell}\ \leq\ \bigl[e^{h/2}%
f_{b_{n}}(e^{-h/c_{b_{n}}})\bigr]^{\ell}.
\end{equation}
As we argued in the derivation of (\ref{FinalBound}), this gives
\begin{equation}
\mathbf{P}(Z_{b_{n}}\leq k_{n}\,|\,Z_{0}=\ell)\,\leq\,e^{-\delta\ell}.
\end{equation}
Consequently, by the Markov property,
\begin{equation}
\mathbf{P}(Z_{n}\leq k_{n})\,\leq\,f_{n-b_{n}}(e^{-\delta}).
\end{equation}
Taking logarithm and using (\ref{not.B}), we obtain (\ref{UB}).

Let us verify the lower bounds in Theorem~\ref{T.Boettcher}. By
(\ref{BoettcherTP}),
\begin{equation}
\mathbf{P}(Z_{n}=k_{n})\ \geq\ \mathbf{P}(Z_{n-b_{n}}=\mu^{n-b_{n}%
})\,\mathbf{P}(Z_{b_{n}}=k_{n}\,|\,Z_{0}=\mu^{n-b_{n}}).
\label{BoettcherLowBound}%
\end{equation}
{F}rom (\ref{CramTrans}),
\begin{equation}
\mathbf{P}(Z_{b_{n}}=k_{n}\,|\,Z_{0}=\mu^{n-b_{n}})\ >\ \bigl[f_{b_{n}%
}(e^{-h/c_{b_{n}}})\bigr]^{\ell_{n}}\,\mathbf{P}\!\left(  _{\!_{\!_{\,}}%
}S_{\ell_{n}}(h,b_{n})=k_{n}\right)  \!, \label{CramTrans'}%
\end{equation}
where $\ell_{n}=\mu^{n-b_{n}}$.

Consider the equation
\begin{equation}
c_{b_{n}}^{-1}\mathbf{E}X_{1}(h,b_{n})\ =\ \frac{f_{b_{n}}^{\prime
}(e^{-h/c_{b_{n}}})\,e^{-h/c_{b_{n}}}}{c_{b_{n}}f_{b_{n}}(e^{-h/c_{b_{n}}}%
)}\ =\ x. \label{LD-equation}%
\end{equation}
Evidently,
\begin{equation}
\frac{f_{b_{n}}^{\prime}(e^{-h/c_{b_{n}}})\,e^{-h/c_{b_{n}}}}{f_{b_{n}%
}(e^{-h/c_{b_{n}}})}\,\Big|_{h=0}\ =\ m^{b_{n}}%
\end{equation}
and
\begin{equation}
\frac{f_{b_{n}}^{\prime}(e^{-h/c_{b_{n}}})\,e^{-h/c_{b_{n}}}}{f_{b_{n}%
}(e^{-h/c_{b_{n}}})}\,\Big|_{h=\infty}\ =\ \mu^{b_{n}}.
\end{equation}
{F}rom these identities and monotonicity of $f_{b_{n}}^{\prime}(e^{-h/c_{b_{n}%
}})\,e^{-h/c_{b_{n}}}/f_{b_{n}}(e^{-h/c_{b_{n}}})$ it follows that
(\ref{LD-equation}) has a unique solution $h_{n}(x)$ for $\mu^{b_{n}}c_{b_{n}%
}^{-1}<x<m^{b_{n}}c_{b_{n}}^{-1}$. Analogously one shows that the equation
$\varphi^{\prime}(h)/\varphi(h)=-x$ has also a unique solution $h(x)$. By the
integral limit theorem (\ref{ILT}), the right-hand side in (\ref{LD-equation})
converges to $-\varphi^{\prime}(h)/\varphi(h)$ and consequently,
$h_{n}(x)\rightarrow h(x)$ as $n\uparrow\infty$. Further, by (\ref{a_nbounds}%
),
\begin{equation}
\frac{\mu}{2m}\ \leq\ x_{n}\ :=\ \frac{k_{n}}{c_{b_{n}}\ell_{n}}\ \leq
\ \frac{1}{2}\,.
\end{equation}
Thus,
\begin{equation}
h(\mu/2m)\ \leq\ \liminf_{n\uparrow\infty}h_{n}\ \leq\ \liminf_{n\uparrow
\infty}h_{n}\ \leq\ h(1/2),
\end{equation}
where $h_{n}:=h_{n}(x_{n})$. It means that there exist $h_{\ast}$ and
$h^{\ast}$ such that $h_{\ast}\leq h_{n}\leq h^{\ast}$ for all $n\geq1$. From
the definition of $\,h_{n}$\thinspace\ and (\ref{LD-equation}) immediately
follows that $\,\mathbf{E}S_{\ell_{n}}(h_{n},b_{n})=k_{n\,}.$\thinspace\ Thus,
applying Lemma~\ref{LCLT}, we get%
\begin{equation}
\lim_{n\uparrow\infty}\bigg|\ell_{n}^{1/2}\,\sigma(h_{n},b_{n})\,\mathbf{P}%
\!\left(  _{\!_{\!_{\,}}}S_{\ell_{n}}(h_{n},b_{n})=k_{n}\right)  \,-\,\frac
{d}{\sqrt{2\pi}}\bigg|\ =\ 0.
\end{equation}
Recall that by (\ref{VarBounds}) we have $\,\sigma(h_{n},b_{n})\geq
C$\thinspace$c_{b_{n\,}}.$\thinspace\ Hence,
\begin{equation}
\liminf_{n\uparrow\infty}\ell_{n}^{1/2}c_{b_{n}}\,\mathbf{P}\!\left(
_{\!_{\!_{\,}}}S_{\ell_{n}}(h_{n},b_{n})=k_{n}\right)  \ \geq C\ >\ 0.
\end{equation}
Moreover, since $\,f_{b_{n}}(e^{-h_{n}/c_{b_{n}}})\geq f_{b_{n}}(e^{-h^{\ast
}/c_{b_{n}}})$\thinspace\ and $\,f_{j}(e^{-h^{\ast}/c_{j}})\rightarrow
\mathbf{E}e^{-h^{\ast}W}>0,$\thinspace\ there exists a $\,\theta>0$%
\thinspace\ such that
\begin{equation}
f_{b_{n}}(e^{-h/c_{b_{n}}})\ \geq\ \theta\label{f-bound}%
\end{equation}
for all $\,n.$\thinspace\ Applying these bounds to the right-hand side in
(\ref{CramTrans'}), we find that
\begin{equation}
\liminf_{n\uparrow\infty}\mu^{b_{n}-n}\log\left[  c_{n}\,\,\mathbf{P}%
(Z_{b_{n}}=k_{n}\,|\,Z_{0}=\mu^{n-b_{n}})\right]  \ \geq\ -C.
\end{equation}
Using this inequality and (\ref{not.B}) to bound the right-hand side in
(\ref{BoettcherLowBound}), we conclude that
\begin{equation}
\liminf_{n\uparrow\infty}\mu^{b_{n}-n}\log\!\left[  _{\!_{\!_{\,}}}%
c_{n}\,\mathbf{P}(Z_{n}=k_{n})\right]  \ \geq\ -C,
\end{equation}
i.e. (\ref{LB}) is proved.

Next we want to extend this result to $\,\mathbf{P}(Z_{n}\leq k_{n}%
).$\thinspace\ Obviously,
\begin{equation}
\mathbf{P}(Z_{n}\leq k_{n})\ \geq\ \mathbf{P}(Z_{n-b_{n}}=\ell_{n}%
)\,\mathbf{P}(Z_{b_{n}}\leq k_{n}\,|\,Z_{0}=\ell_{n}).
\end{equation}
Then, using (\ref{CramTrans}) with $h=h_{n\,},$ we have
\begin{equation}
\mathbf{P}(Z_{n}\leq k_{n})\ \geq\ \mathbf{P}(Z_{n-b_{n}}=\ell_{n}%
)\,\bigl[f_{n}(e^{-h_{n}/c_{b_{n}}})\bigr]^{\ell_{n}}\,\mathbf{P}\!\left(
_{\!_{\!_{\,}}}S_{\ell_{n}}(h,b_{n})\leq k_{n}\right)  \!.
\end{equation}
By the central limit theorem,
\begin{equation}
\lim_{n\uparrow\infty}\mathbf{P}\!\left(  _{\!_{\!_{\,}}}S_{\ell_{n}}%
(h,b_{n})\leq k_{n}\right)  \ =\ \frac{1}{2}\ .
\end{equation}
{F}rom this statement and (\ref{f-bound}) it follows that
\begin{equation}
\liminf_{n\uparrow\infty}\mu^{b_{n}-n}\log\mathbf{P}(Z_{n}\leq k_{n}%
)\ \geq\ \liminf_{n\uparrow\infty}\mu^{b_{n}-n}\log\mathbf{P}(Z_{n-b_{n}}%
=\mu^{n-b_{n}})+\log\theta.
\end{equation}
Recalling (\ref{mu^n-probability}), the proof of Theorem~\ref{T.Boettcher} is
complete.\hfill$\square$

\begin{remark}
[\textbf{To the proof of Remark}~\ref{R.w.at.0}]\label{R.w.at.0.details}To
prove (\ref{w.at.0}) one can use the methods from the proof of
Theorem~\ref{T.Boettcher}. But some changes are needed, since in
Remark~\ref{R.w.at.0} we deal with absolutely continuous distributions.

Instead of (\ref{BoettcherTP}) we shall use (\ref{SelfSimilarity}). Putting
there $x=ym^{k}$ and $\,k=k_{y}=\max\{j\geq1:\,m^{j}\leq\mu^{j}/2y\}$%
\thinspace\ we obtain
\begin{equation}
w(y)\ =\ m^{k_{y}}\sum_{\ell=\mu^{k_{y}}}^{\infty}\mathbf{P}(Z_{k_{y}}%
=\ell)\,w^{\ast\ell}(ym^{k_{y}}). \label{SelfSimilarity1}%
\end{equation}
For every $h\geq0$ we may define the density function
\begin{equation}
w_{h}(x)\ :=\ \frac{e^{-hx}}{\varphi(h)}\ w(x), \label{2R}%
\end{equation}
corresponding to the Cram\'{e}r transform of $\,W$. By Lemma~\ref{w(x)Bound},
$\,C_{w}:=\sup_{x\geq0}w(x)<\infty$\thinspace\ in the present B\"{o}ttcher
case. Hence, $\,\sup_{x\geq0}w_{h}(x)\leq C_{w}/\varphi(h).$\thinspace\ By
induction (analogously to Lemma~\ref{ConFunction}),%
\begin{equation}
\sup_{x\geq0}w_{h}^{\ast\ell}(x)\ \leq\ \frac{C_{w}}{\varphi(h)}\,,\qquad
\ell\geq1. \label{SupNorm}%
\end{equation}
It is easy to see that
\begin{equation}
w_{h}^{\ast l}(x)\ =\ \frac{e^{-hx}}{\varphi^{\ell}(h)}\ w^{\ast\ell
}(x),\qquad\ell\geq1. \label{DensityCramTrans}%
\end{equation}
From this identity and (\ref{SupNorm}) it follows that
\begin{equation}
w^{\ast\ell}(x)\,\leq\,C_{w}\,\varphi^{\ell-1}(h)\,e^{hx}.
\end{equation}
Therefore, for all $\,\ell\geq\mu^{k_{y}},$
\begin{equation}
w^{\ast\ell}(ym^{k_{y}})\ \leq\ \frac{C_{w}}{\varphi(h)}\ \Bigl[e^{hym^{k_{y}%
}/\mu^{k_{y}}}\varphi(h)\Bigr]^{\ell}.
\end{equation}
Further, by the definition of $\,k_{y\,}$,
\begin{equation}
\frac{\mu}{2my}\,\leq\,\frac{m^{k_{y}}}{\mu^{k_{y}}}\,\leq\,\frac{1}{2y}\,,
\label{k_y-bounds}%
\end{equation}
and consequently,
\begin{equation}
w^{\ast\ell}(ym^{k_{y}})\ \leq\ \frac{C_{w}}{\varphi(h)}\ \bigl[e^{h/2}%
\varphi(h)\bigr]^{\ell}.
\end{equation}
Before (\ref{FinalBound}) we showed that $\,e^{h_{0}/2}\varphi(h_{0})\leq
e^{-\delta}.$\thinspace\ As a result we have the bound
\begin{equation}
w^{\ast\ell}(ym^{k_{y}})\ \leq\ \frac{C_{w}}{\varphi(h_{0})}\ e^{-\delta\ell}.
\end{equation}
Entering this into (\ref{SelfSimilarity1}) gives
\begin{equation}
w(y)\,\leq\,C\,m^{k_{y}}\,f_{k_{y}}(e^{-\delta}).
\end{equation}
Taking logarithm and using (\ref{not.B}), we see that
\begin{equation}
\limsup_{y\rightarrow0}\,\mu^{-k_{y}}\log w(y)\ \leq\ \log B(e^{-\delta}).
\label{UpperBound}%
\end{equation}

Now we deal with a corresponding lower bound of $\log w(y)$. By
(\ref{SelfSimilarity1}) and (\ref{DensityCramTrans}),
\begin{align}
w(y)  &  >m^{k_{y}}\,\mathbf{P}(Z_{k_{y}}=\mu^{k_{y}})\,w^{\ast\mu^{k_{y}}%
}(ym^{k_{y}})\nonumber\label{LowBound}\\
&  >\mathbf{P}(Z_{k_{y}}=\mu^{k_{y}})\,\varphi^{\mu^{k_{y}}}(h)\,w_{h}%
^{\ast\mu^{k_{y}}}(ym^{k_{y}}),\qquad h>0.
\end{align}
Recalling that $h(x)$ is the unique solution of the equation $\,\varphi
^{\prime}(h)/\varphi(h)=-x$\thinspace\ and using (\ref{k_y-bounds}), one gets
the inequality $\,h(ym^{k_{y}}/\mu^{k_{y}})\leq h(\mu/2m).$\thinspace\ Thus,
by monotonicity of $\,\varphi$,
\begin{equation}
\varphi^{\mu^{k_{y}}}\left(  _{\!_{\!_{\,}}}h(ym^{k_{y}}/\mu^{k_{y}})\right)
>\,\varphi^{\mu^{k_{y}}}\left(  _{\!_{\!_{\,}}}h(\mu/2m)\right)  =\,\exp
[-C\mu^{k_{y}}]. \label{phi-bound}%
\end{equation}
If in (\ref{2R}) we set $\,h=h(ym^{k_{y}}/\mu^{k_{y}}),$\thinspace\ then
$\,w_{h}^{\ast\mu^{k_{y}}}(ym^{k_{y}})$\thinspace\ is the value of the density
function of the sum $\,\sum_{j=1}^{\mu^{k_{y}}}W_{j}(-h)$\thinspace\ at the
point $\,\mathbf{E}\sum_{j=1}^{\mu^{k_{y}}}W_{j}(-h).$\thinspace\ Thus, by the
central limit theorem for densities (\cite[Theorem~VII.7]{Petrov1975}),
\begin{equation}
\lim_{y\rightarrow0}w_{h}^{\ast\mu^{k_{y}}}(ym^{k_{y}})\ =\ \frac{1}%
{\sqrt{2\pi}}\,. \label{DensityCLT}%
\end{equation}
Putting $h=h(ym^{k_{y}}/\mu^{k_{y}})$ in (\ref{LowBound}) and using
(\ref{mu^n-probability}), (\ref{phi-bound}), and (\ref{DensityCLT}), we
obtain
\begin{equation}
\liminf_{y\rightarrow0}\,\mu^{-k_{y}}\log w(y)\ \geq\ -C. \label{LowerBound}%
\end{equation}
Combining (\ref{UpperBound}) and (\ref{LowerBound}) we get
\begin{equation}
\log w(y)\,\asymp\,-\mu^{k_{y}}.
\end{equation}
Then the relation $\mu^{k_{y}}\asymp y^{-\beta/(1-\beta)}$ finishes the
proof.\hfill$\Diamond$
\end{remark}

{\small
\bibliographystyle{plain}
\bibliography{bibtex,bibtexmy}
}%

%

\end{document}